\documentclass[12pt]{amsart}
\usepackage{latexsym,fancyhdr,amssymb,color,amsmath,amsthm,graphicx,listings,comment}
\usepackage[section]{placeins}
\definecolor{red1}{rgb}{1,0.9,0.9} \definecolor{blue1}{rgb}{0.9,0.9,1} \definecolor{green1}{rgb}{0.9,1,0.9} 
\definecolor{yellow1}{rgb}{1,1,0.6} \definecolor{yellow2}{rgb}{1,1,0.8}
\definecolor{brightgreen}{rgb}{0.8,1,0.8}
\definecolor{brightblue}{rgb}{0.9,0.9,1}
\definecolor{brightyellow}{rgb}{1.0,1.0,0.6}
\def\ssection#1{ \begin{center} \vspace{1mm} \fcolorbox{brightgreen}{brightgreen}{ \parbox{14.2cm}{\bf{ \large{ #1 }}}} \vspace{1mm} \end{center} }
\def\ttheorem#1{ \begin{center} \vspace{0mm} \fcolorbox{brightyellow}{brightyellow}{ \parbox{14.2cm}{\bf{ \large{ #1 }}}} \vspace{0mm} \end{center} }
\def\question#1{ \vspace{2mm} \begin{center} \fcolorbox{green1}{green1}{ \parbox{14.2cm}{{\bf Question:} #1}} \vspace{2mm} \end{center} }

\def\remark#1{ \vspace{2mm} \begin{center} \fcolorbox{yellow1}{yellow1}{ \parbox{14.2cm}{{\bf Remark:} #1}} \vspace{2mm} \end{center} }

\let\paragraph\subsection

\title{The Babylonian graph}
\author{Oliver Knill}
\date{June 25, 2022}
\address{Department of Mathematics \\ Harvard University \\ Cambridge, MA, 02138 }

\begin{document}
\maketitle

\begin{abstract}
The Babylonian graph $B$ has the positive integers as vertices and connects
two if they define a Pythagorean triple. Triangular subgraphs correspond to Euler bricks. 
What are the properties of this graph? Are there tetrahedral subgraphs corresponding to Euler tesseracts?
Is there only one infinite connected component? Are there two Euler bricks in the graph that are disconnected?
Do the number of edges or triangles in the subgraph generated by the first n vertices grow 
like of the order n W(n), where n is the product log? We prove here some simple results.
In an appendix, we include handout from a talk on Euler cuboids given in the year 2009.
\end{abstract}

\section{Babylonian graphs}

\begin{figure}[!htpb]
\scalebox{1.0}{\includegraphics{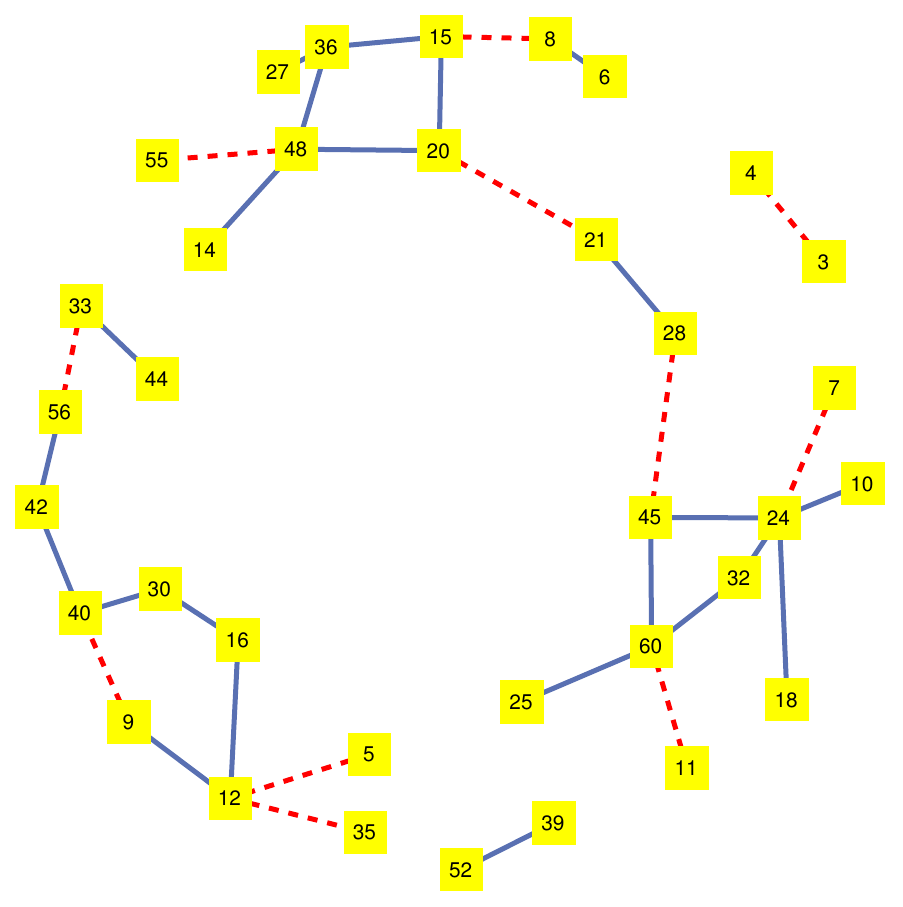}}
\label{adc}
\caption{
The graph $B_{60}$ without isolated vertices like $\{1\}$ or $\{2\}$. 
The edges corresponding to primitive Pythagorean triples are dashed. 
}
\end{figure}

\paragraph{}
For every positive integer $n$, let {\bf Babylon-n} denote the simple graph $B_n=(V_n,E_n)$
with vertex set $V_n=\{1,\dots,n\}$ and edge set $E_n=\{ (a,b), a^2+b^2 \in \mathbb{N},a,b \leq n \}$.
We have a nested increasing sequence of graphs $B(n)$ which starts with $B_0=K_1$ and leads to
the {\bf Babylonian graph} $B=(\mathbb{Z}^+,E)$ with {\bf Pythagorean triples} $(a,b)$ as edges 
is infinite. Triangles, complete subgraphs $K_3$ in $B$ correspond to {\bf Euler bricks} 
\cite{dicksonII} (chapter XIX, see also \cite{KnillTreasure}. 
The more extended notes \cite{KnillEndpaper} is attached here as an Appendix).
A subclass of Euler bricks can be obtained by parametrizations like the {\bf Saunderson parametrization} 
$a=u (4v^2-w^2), b=v (4u^2-w^2), c=4 u v w$. Triangles $(a,b,c)$ for which additionally the sum $a^2+b^2+c^2$ 
is an integer correspond to {\bf perfect Euler bricks}, an object which has not yet been found and which might
not exist. The graphs $B_n$ are the from $V_n=\{1,\dots,n\}$ induced subgraph of $B$. 
Each $B_n$ is a subgraph of $B_{n+1}$ and the limit $B=(\mathbb{Z}^+,E)$ encodes all Pythagorean triples.

\begin{figure}[!htpb]
\scalebox{1.0}{\includegraphics{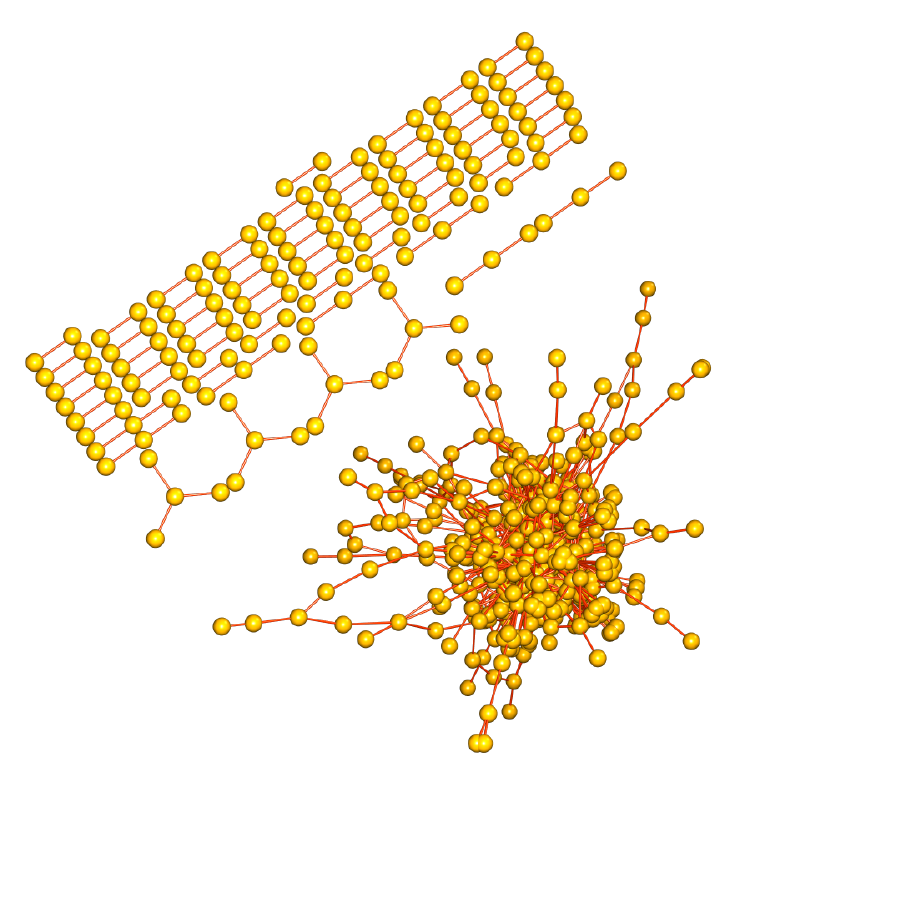}}
\label{adc}
\caption{
The graph $B_{1000}$ without the $0$-dimensional isolated points. 
}
\end{figure}

\paragraph{}
We are interested in the {\bf largest connected component} $B_n'$ of $B_n$ and
in the connectivity or symmetry properties of $B$. Are there $K_4$ subgraphs in $B$?
What groups act as graph isomorphisms on $B$?
While $B$ has some small components like the single vertex $\{1\}$ or the 
single vertex $\{2\}$ or the isolated $K_2$ graph $\{3 \leftarrow 4 \}$, we expect
that $B$ only has one large infinite connectivity component.

\begin{figure}[!htpb]
\scalebox{0.7}{\includegraphics{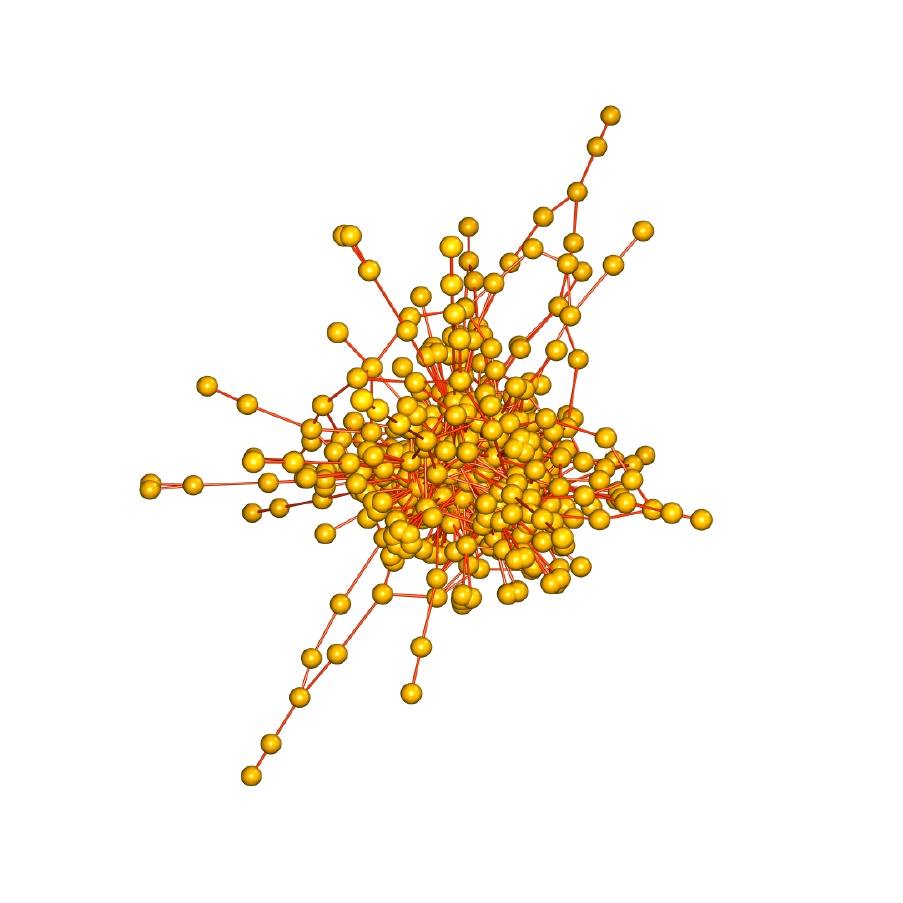}}
\label{babylon1000}
\caption{
The main connected component $B_{1000}'$ 
of the Babylonian graph $B_{1000}$. Its $f$-vector is
$f=(480, 952, 10)$, its Euler characteristic is $\chi(B_{1000}')=480-952+10=-462$
and its Betti vector is $b=(b_0,b_1)=(1,463)$. The
graph $B_{10000}$ itself has $f$-vector $f=(1000, 1034, 10)$, Betti vector
$b=(b_0,b_1) = (439, 463)$ and $\chi(B_{10000})=f_0-f_1+f_2=b_0-b_1=-24$. 
}
\end{figure}

\paragraph{}
The question of existence of perfect Euler bricks appears to be difficult. The popularity of the
problem persists. It is now also in a list of problems in \cite{VisionsOfInfinity}
discussed beyond the Millenium problems. Euler bricks are triangles in $B$ and
correspond to points $(x,y,z)$ in $\mathbb{R}^3$ located on three 
cylinders $x^2+y^2=a^2,y^2+z^2=b^2,z^2+x^2=c^2$ 
with integer radii $a,b,c$. In an Euler brick triangle, it is not possible 
that all three pairs are primitive. The graph 
induced by the primitive edges has no triangles.  Figure~\label{perfect} illustrates the geometric problem 
to find integer points $(x,y,z)$
on the intersection of three perpendicular main axes-centered cylinders with integer radius. The perfect Euler
brick problem is to find such points which also have integer distance to the origin. 
Also the problem of {\bf Euler tesseracts} can be seen geometrically. It is 
the problem to find the intersection of six perpendicular $3D$-cylinders 
$x_i^2+x_j^2=r_{ij}^2$, $1 \leq i<j \leq 4$ with 
integer radius $r_{ij}$ in $\mathbb{R}^4$. It rephrases to find $K_4$ subgraphs of $B$. 
Whether this is possible is not clear. 

\paragraph{}
There are quite many Diophantine problems using squares. One similar looking problem is the 
{\bf Mengoli six square problem} of Mietro Mengoli who by the way also would suggest the {\bf Basel problem}, to find the value of 
$\sum_{k=1}^{\infty} 1/k^2$. Mengoli asked for triples $x \leq y \leq z$ of integers such that the sum and difference of any two are
squares. In other words,  $x+y,y+z,x+z,y-x,z-y,z-x$ should all be squares. Euler found the smallest one.
The solution found by Ozanam in 1691 is $(1873432,2288168,2399057)$. \cite{NastasiScimone}

\begin{figure}[!htpb]
\scalebox{0.5}{\includegraphics{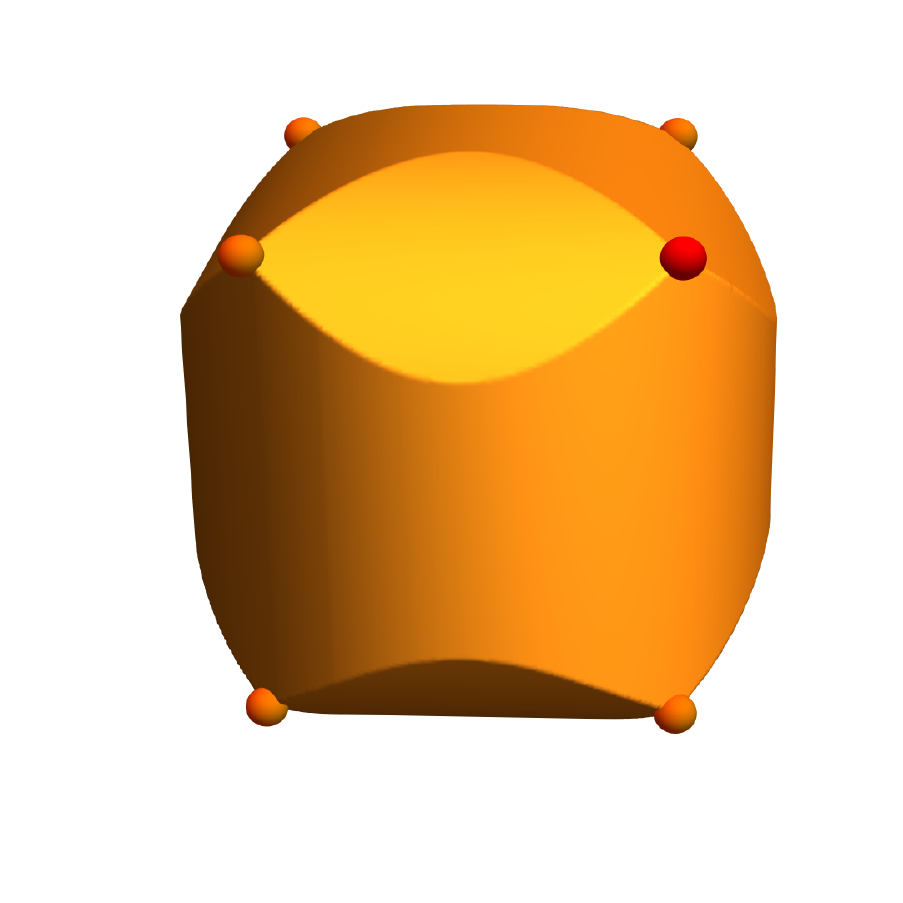}}
\scalebox{0.5}{\includegraphics{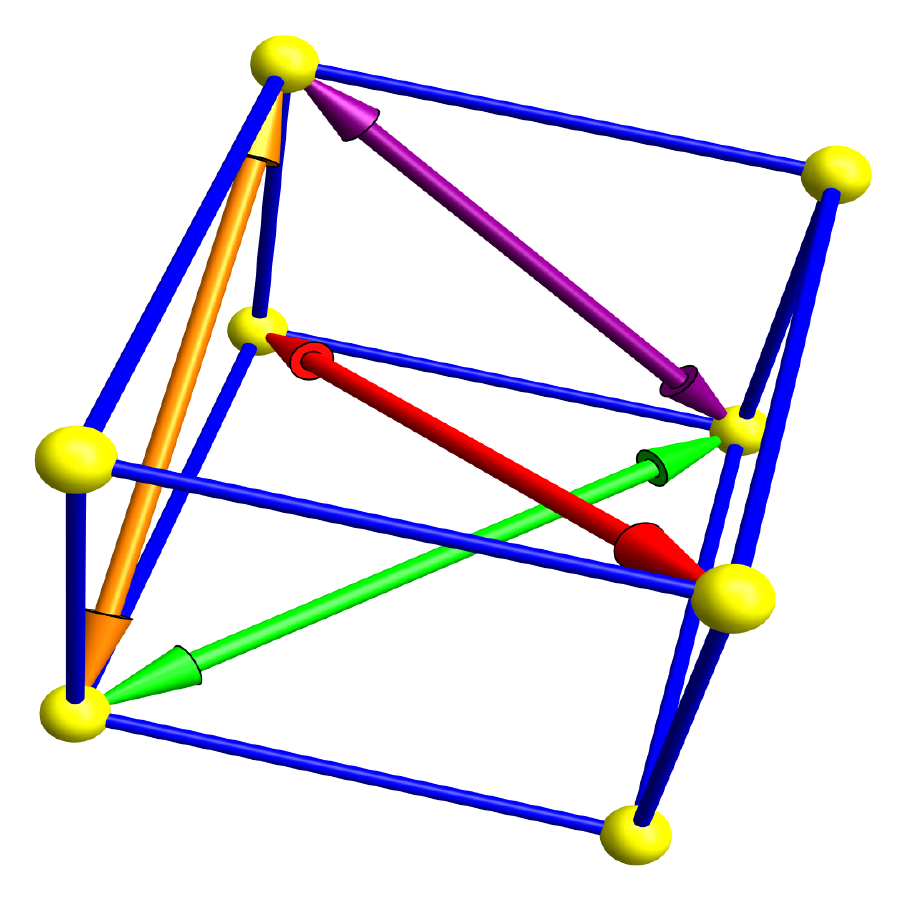}}
\label{perfect}
\caption{
The intersection of three cylinders with integer radius defines Euler bricks.
If the intersection point $(x,y,z)$ also has integer distance to the origin, 
we have a perfect Euler brick. To the right we see an Euler brick. 
}
\end{figure}

\section{Four Questions}

\paragraph{}
Here are some natural questions. They can also be formulated as conjectures. Despite having
this posted since February 2022, it is probably safer to still keep it as questions and not upgrade it 
to conjectures.
One reason is that one or the other question could turn out to be obvious. When studying a problem for
the first time, it is possible to miss something obvious. 
It can be that one or the other question are already answered in an other context or is 
a special case of a general theorem. We had looked a couple of years ago at the literature of Euler bricks 
while preparing for a math circle talk. That handout in a treasure hunting theme \cite{TreasureIsland}
is here attached in an appendix. 

\paragraph{}
First of all, we believe that there is only one largest connected component $B_n'$ for all $n$ and 
that it will end up to be a {\bf single infinite component} in the Babylonian graph $B$. 
In principle, it is not yet excluded that there are several infinite disconnected components of $B$.
Question \textcircled{A} asks whether such an ``eternal maximal component" exists.

\question{ \textcircled{A} Does $B$ have only one infinite connected component?}

\paragraph{}
Question \textcircled{B} is about the existence of 
{\bf Hyper Euler bricks} or {\bf Euler tesseracts}. It addresses the question about the 
{\bf maximal dimension} of $B$. While unlikely, it would in principle be possible that the 
maximal dimension is infinity, meaning that there are complete subgraphs $K_n$ of $B$ for
any integer $n$. Already the question of three dimensional complete subgraphs $K_4$ is unclear:

\question{ \textcircled{B} Is there a $K_4$ sub-graph in $B$?  }

{\bf Euler tesseracts} are hyper cubes $\{ (x_1,x_2,x_3,x_4) \in \mathbb{R}^4, 0 \leq x_1 \leq x,
0 \leq x_2 \leq y, 0 \leq x_3 \leq z, 0 \leq x_4 \leq w \}$ with integer $x,y,z,w$ side length
for which all the 6 2D-face diagonals have integer length. 
One of the ways to show that a system of Diophantine equations has no solution is to look 
at the system modulo a prime $p$. If there is no solution modulo $p$, then there is no solution 
at all. In the {\bf tesseract problem} we have to solve the system of Diophantine equations
$$ x^2+y^2=a^2, y^2+z^2=b^2, z^2+w^2=c^2, w^2+x^2=d^2, x^2+z^2=e^2,y^2+w^2=f^2 \; $$
for the 10 integer variables $x,y,z,w,a,b,c,d,e,f \in \mathbb{Z}^+ =\{1,2,3, \dots \}$. 
We have not started to look for solutions yet, but the strategy is similar than when looking
for perfect Euler bricks: take a parametrization $(x,y,z)$ of Euler bricks like the Saunderson 
parametrization, then we have a function $F_{x,y,z}(w) = d(\sqrt{x^2+w^2}) + d(\sqrt{y^2+w^2}) + d(\sqrt{z^2+w^2})$
where $d(t)$ is the distance to the nearest integer. Now, for large $x,y,z$ and $w$, 
the dynamics $F_{x,y,z}(w) \to F_{x,y,z}(w+1)$ is by linear approximation close to a
translation $t \to t+\alpha$ which then by a continued fraction expansion allows to find
$w$ for which $F_{x,y,z}(w)$ is very small and since the possible distances are quantized, once we
are close enough, we hit a solution. This is how we have searched for perfect Euler bricks. 
A more sophisticated search using a multi-dimensional approach, leading to multi-variable
Chinese remainder type problems \cite{Knill2012}.

\paragraph{}
Question \textcircled{C} asks about the maximal dimension of the non-major connected components. 
Example sub-graphs $\{1\},\{2\},\{3,4\}$ remain separated also in $B$. We have not yet seen
any example, where a non-major connected component has maximal dimension larger than $1$: 

\question{ \textcircled{C} Are there two connected components with triangles? }

In other words, we look for pairs of Euler bricks, so that there is no
connection from one brick to the other brick in $B$. We have also not seen any in any $B_n$. 
It could be possible that there exists a $B_n$ with two disconnected $K_3$ subgraphs which 
however get reunited in a larger $B_m$ for $m >n$. .

\paragraph{}
The fourth and last major question \textcircled{D} is a particular {\bf growth rate question}. 
Let $W(x)$ denote the {\bf product log} = {\bf Lambert W-function} which is defined as the inverse of the
function $y=x \log(x)$.  This function naturally occurs in the {\bf prime number theorem} which tells that the
n'th prime $p_n$ is of the order $p_n \sim n \log(n)$ meaning $W(p_n) \sim n$. Now the 
{\bf number of primitive edges} in $B_n$ grow like $C_1 n$, where $C_1$ is a concrete number 
expressible as an area of the parameters lattice points $(u,v)$ in 
$[0,\sqrt{n}] \times [0,\sqrt{n}]$ such that $\Phi(u,v) = (u^2-v^2,2uv) \in [0,n] \times [0,n]$. 
Having a growth $C_1 n$ of  primitive edges, we expect $C_1 n W(n)$ to be the growth of all edges. 

\question{ \textcircled{D}  Does $f_1(B_n)/(n W(n))$ converge? }

\paragraph{}
There are many other quantities one could look at.
We can look at the {\bf f-vector} $(f_0,f_1,f_2,\dots) = (n,f_1(B_n),f_2(B_n),\dots)$ or the Betti numbers 
$b_k(B_n) = {\rm dim}({\rm ker}(L_k(B_n))$ where $L_k(B_n)$ is the {\bf $k$-form Laplacian}.
Of interest are the number of connected components $b_0(B_n)$, the number $b_1(B_n)$, a genus, which 
measures of the number of one dimensional ``holes" in $B_n$,
the maximal vertex degree, the distribution of the vertex degrees, the growth of the Euler characteristic 
$\chi(B_n)=f_0(B_n)-f_1(B_n)+f_2(B_n)- \dots = b_0(B_n)-b_1(B_n)+b_2(B_n)- \dots$, the inductive dimension of 
$B_n$ or $B_n'$ with the ultimate goal to give lower and upper bounds. 
We looked first numerically at the diameter change up to $n=10000$. 
The diameter goes to infinity because ${\rm Diam}(B_1(5000))=18$, ${\rm Diam}(B_1(10000))=29$.
A logarithmic growth which is justified by scaled graphs 
$x_1,x_2, \dots, x_n = M x_1, M x_2, M x_3, \dots, M x_n$ etc.

\begin{figure}[!htpb]
\scalebox{0.7}{\includegraphics{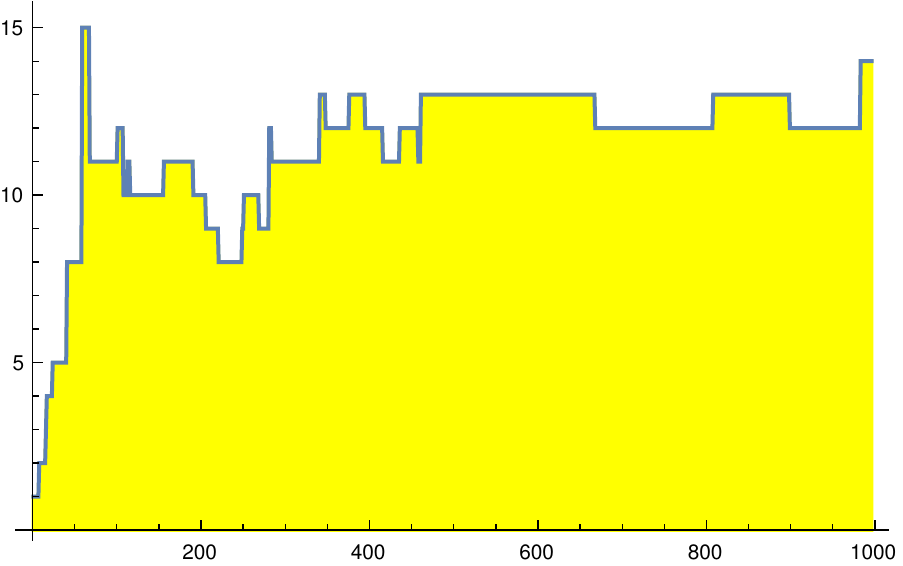}}
\label{diameter}
\caption{
The diameter of $B_n$ up to $n=1000$.
}
\end{figure}

\paragraph{}
Whenever we have a graph $G$, the {\bf graph complement} $\overline{G}$ is of interest.
The operation of taking graph complements
is an {\bf involution} on the class of all graphs.
Even for very simple graphs, the graph complement can be interesting.
See \cite{GraphComplements} for cyclic or linear graphs, where graph complements are either contractible
or homotopic to spheres or wedge sums of spheres. What are the properties of the graph complement of $B$?

\section{Low hanging fruits}

\paragraph{}
One can wonder first whether there are {\bf isolated vertices} that are not connected to anything else. 
These are zero dimensional connected components of the graph. 
There are exactly two vertices with this property. We know that $1+b^2, 4+b^2$ are never 
a square. 

\remark{ \textcircled{1} There are exactly two isolated single vertices $1,2$ in $B$.}

{\bf Proof:} Already the {\bf primitive Babylonian graph} $B_p$ which connects only points $a,b$
if $a,b,\sqrt{a^2+b^2}$ is a primitive Pythagorean triple has no isolated points except $1,2$: 
the reason is that all odd numbers larger than $1$ are of the form $u^2-v^2$ and all 
even numbers larger than $2$ are of the form $2uv$ for some positive distinct $u,v$. 

\paragraph{}
There are other isolated connected components like $\{ 3,4 \}$ belonging to the primitive tripe $3^2+4^2=5^2$. 
This can not be connected to anything else. If $9+b^2 = c^2$, then $c^2-b^2=9$ which 
is only possible for $b=4,c=5$. The relation $16+b^2=c^2$ is only possible for $b^2=9$. 

\paragraph{}
One can wonder about the asymptotic distribution of the vertex degrees. When looking at the
sequence $B_n$ of graphs, there is an increasing part $C_n \subset B_n$ for which the vertex
degrees do no more change when increasing $n$. The reason is 
that the monotone sequence ${\rm deg}(B_n)(x)$ converges:

\remark{ \textcircled{2} ${\rm deg}_B(x)$ is finite for all $x \in V$. }

{\bf Proof:} every edge is a Pythagorean triple which is a multiple of a primitive
Pythagorean triple and so of the form 
$$  (a,b)=(2 p q c,(p^2-q^2) c)  \; , $$ 
where $p,q,c$ are integers. If we fix an integer like $b$, there are only 
finitely many solutions $(p^2-q^2) c = b$ because both $c$  and $p^2-q^2$ have 
to be recruited from factors of $b$ which is finite. For a fixed factor $r$ of $b$
the Diophantine equation $p^2-q^2=r$ has only finitely many solutions because both 
$p,q$ have to be smaller than $\sqrt{r}$. 

\paragraph{}
One can also wonder how many infinite connected components there are $B$. 
We do not know yet. We can get infinite connected components 
for a primitive $(a,b) = (2uv, u^2-v^2)$ if it is connected to a 
multiple of itself. This can indeed happen and proves

\remark{ \textcircled{3} The diameter of $B$ is infinite.}

{\bf Proof:} it is enough to construct a concrete path in $B$ going to infinity: 
There is a connection from $n=5$ to $n=30$ given by 
$5 \to 12 \to 16 \to 30$. This scales. The path $30 \to 40 \to 96 \to 180$
for example extends the other path so that we have a connection from $5$ to $180$. 
We can continue like that and get an infinite path in $B$. 


\paragraph{}
We could call a connected component in which a scaling $(a,b) \to m (a,b)$ exists
a {\bf component with scale symmetry}. We have just seen that such components are
infinite. Are there components which are not scale invariant? 

\paragraph{}
Since $a=u^2-v^2,b=2uv$ parametrizes all {\bf primitive triples},
there is a constant $C$ such that there are asymptotically $C*n$
primitive edges in the graph $B_n$. This limit can be computed explicitly as
we can draw out the region in the parameter domain leading to triples $\leq n$. 
But then we have also non-primitive ones which come from scaled smaller 
primitive ones.

\begin{figure}[!htpb]
\scalebox{0.8}{\includegraphics{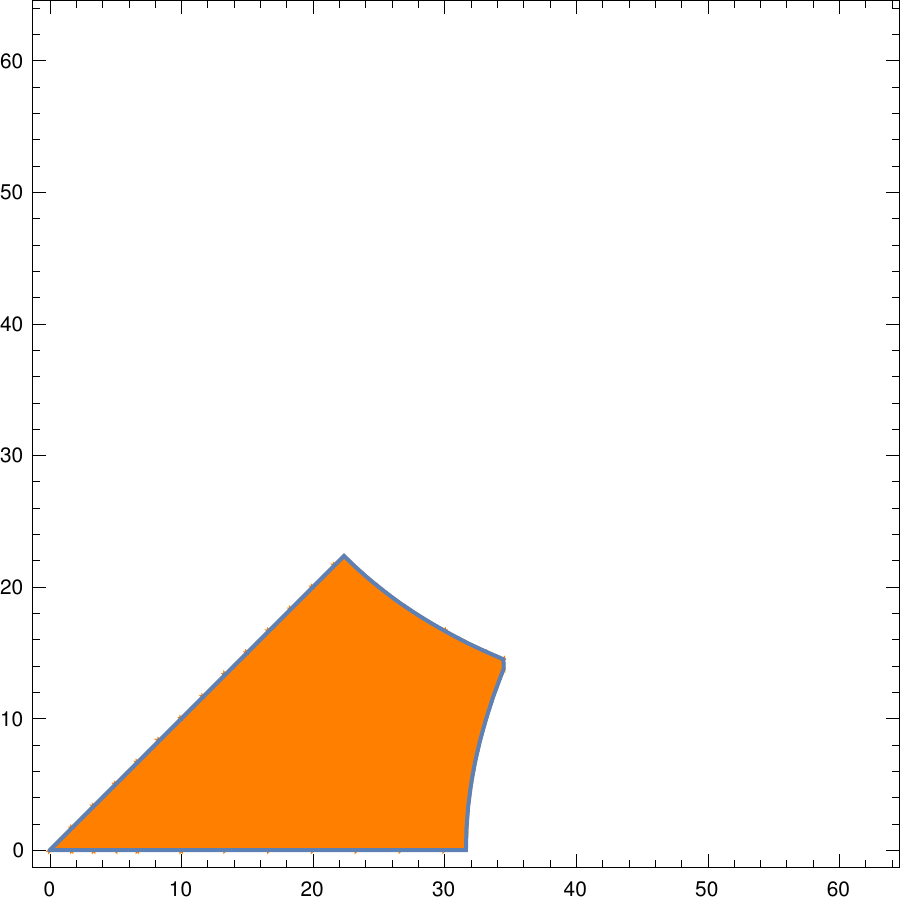}}
\label{uv}
\caption{
The uv-parameter region $- \leq u^2-v^2 \leq n, 0 \leq 2uv \leq n, u \geq 0, v \geq 0$.
For each lattice point $(u,v)$ in that region, we get a
Pythagorean triple $u^2-v^2,2uv,u^2+v^2$ on $[0,2\sqrt{n}] \times [0,2\sqrt{n}]$
and so an edge in $B_n$.
}
\end{figure}

\paragraph{}
Here is an other simple observation about {\bf leafs} in $B$. These are 
vertices $x$ for which the unit sphere $S(x)$ (the subgraph generated by all 
vertices attached to $x$)  contains only one point. 

\remark{ \textcircled{4} If $p$ is an odd prime, then $p$ is a leaf in $B$. }

{\bf Proof:} We must have $p = u^2-v^2$ so that $p$ belongs to the primitive
Pythagorean triple $2uv, u^2-v^2$. Now, $u^2-v^2=p$ implies with $v=u-k$
that $p=u^2-(u-k)^2 = 2uk-k^2$. Since $p$ is prime, $k=1$ meaning that we
have only one choice to solve $u^2-v^2=p$ for $u,v$. The Pythagorean triple
is not fixed. 

\paragraph{}
By definition, the number of vertices $f_0(B_n)=n$. Let $W(x)$ denote the inverse of 
the function $x \to x e^x$. It is called the {\bf Lambert W function}. 
Motivated from the prime number theorem telling that that $n$ grows like $P(n) W(P(n))$, 
where $P(n)$ is the n'th prime number, it is likely that $\lim_{n \to \infty} f_1(B_n)/(n W(n)) = C_1$ 
and $\lim_{n \to \infty} f_2(B_n)/(n W(n)) = C_2$ exist: 

\begin{figure}[!htpb]
\scalebox{0.5}{\includegraphics{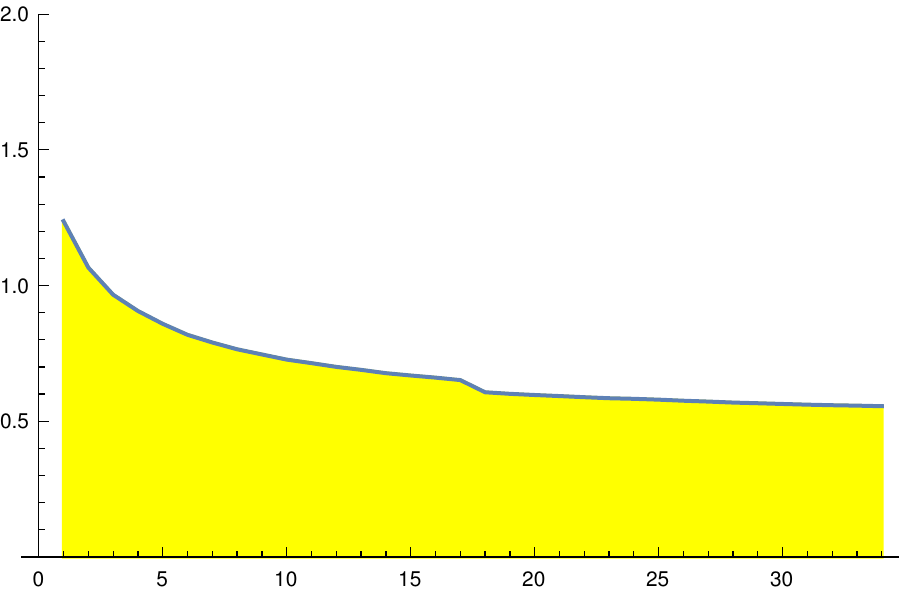}}
\scalebox{0.5}{\includegraphics{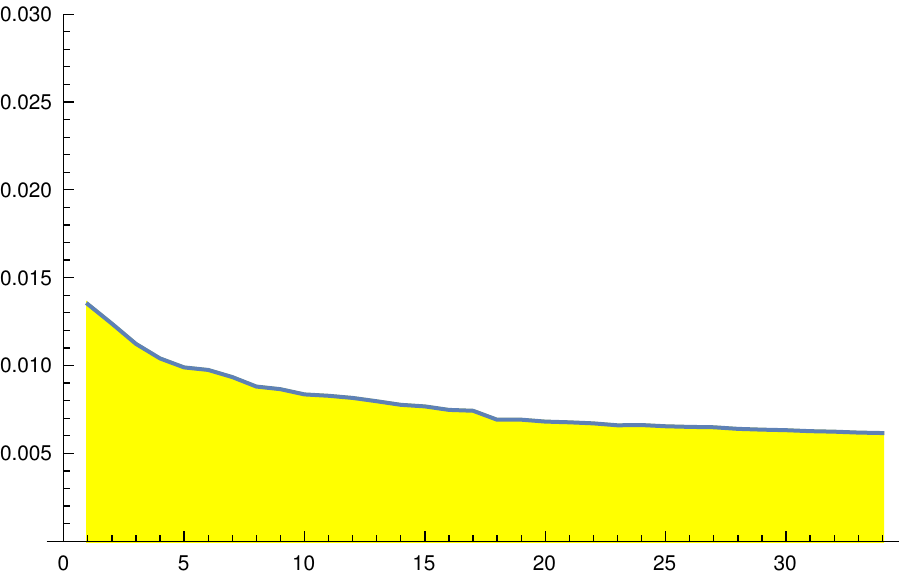}}
\label{f-vector}
\caption{
The number of edges and triangles for the main component.
$ f_1(B_{n*1000}')/(n*W(n)*1000)$ and $f_2(B_{n*1000}')/(n*W(n)*1000)$ 
for $n=1,\dots,25$.  
}
\end{figure}

\begin{figure}[!htpb]
\scalebox{0.5}{\includegraphics{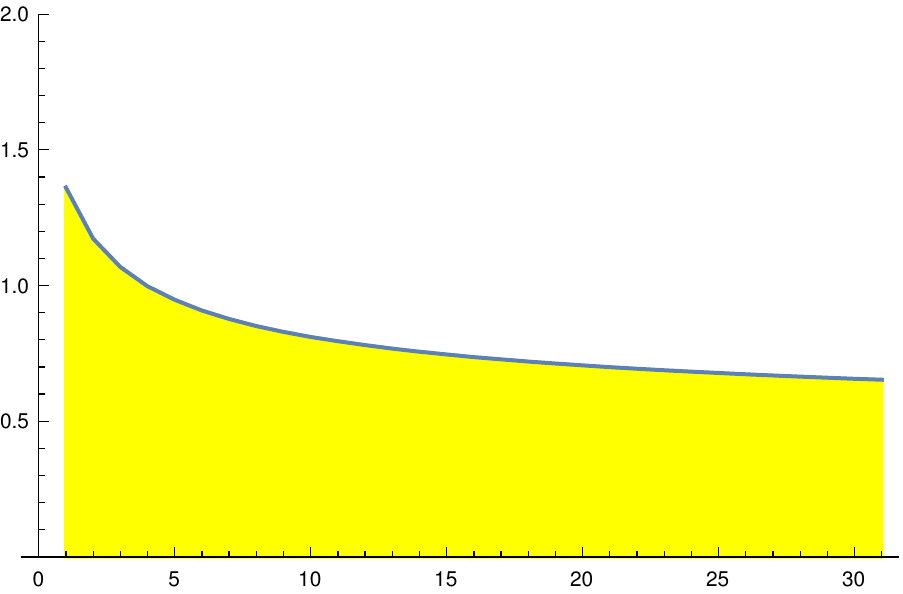}}
\scalebox{0.5}{\includegraphics{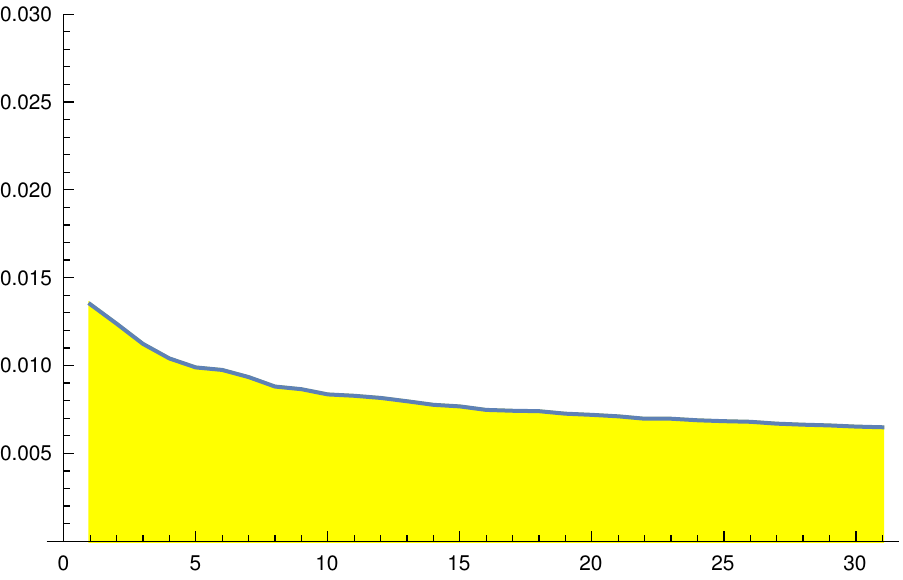}}
\label{f-vector of fullgraph}
\caption{
Number of edges and triangles for the full graph with all components:
$f_1(B_{n*1000})/(n*W(n)*1000)$ and $f_2(B_{n*W(n)*1000})/(n*1000)$ 
for $n=1,\dots ,20$.  
}
\end{figure}

\paragraph{}
This would lead to a result that $C=\lim_{n \to \infty} \chi(B_n)/(n W(n))=C_0-C_1+C_2-C_3+...$ exists. 

\paragraph{}
Here is a result which is somehow interesting. A graph can be defined to be planar if it does not 
contain a homeomorphic copy of $K_5$ or $K_{3,3}$. (There is also the traditional
definition of planar using a topological embedding on a 2-sphere, but the just given one is equivalent by 
{\bf Kuratowski's theorem}. The combinatorial definitio
has has the advantage that is purely graph theoretical and does not refer to topology of Euclidean space.

\remark{ \textcircled{5} $B_n$ is planar if and only if $n \leq 95$}

\begin{proof} 
Since $B_{95}$ is planar all $B_n$ with $n \leq 95$ are planar. 
Since $B_{96}$ is non planar, all $B_n$ with $n \geq 96$ are non-planar.
\end{proof} 

\begin{figure}[!htpb]
\scalebox{0.9}{\includegraphics{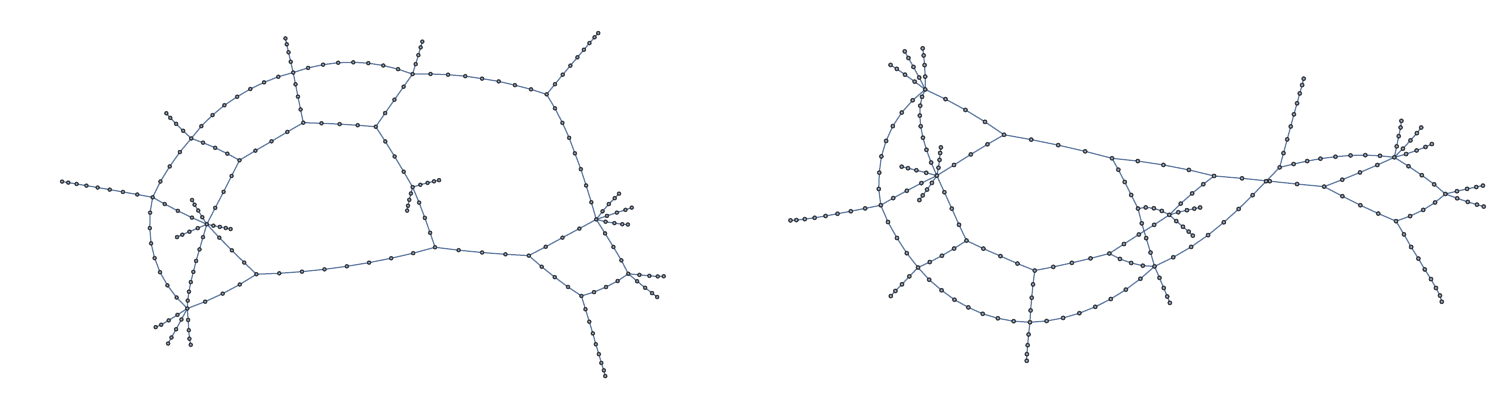}}
\label{B(95) and B(96)}
\caption{
The graph $B_{95}$ (looking like a piglet) is planar, the graph $B_{96}$ (looking like a chicken) is 
not. We show here the second Barycentric refinements of the main connected components. The non-planar property 
of $B_{96}$ is not well visible. We have to look closely at the crossing of edges. There is a crossing 
which is not a node. 
}
\end{figure}

\section{Babylonian triplets and Pythagoras}

\paragraph{}
Pythagorean triples appeared in Babylonians tablets. The most famous is
Plimpton 322 \cite{Robson2002,Mansfield2021}. An other one  is Si.427. We could call the examples of Pythagorean
triples which appear in some Babylonian text a {\bf Babylonian triplet}. 
The major known clay documents which list Babylonian triples appear all during the time 1900-1600 BC.
Pythagorean triples also appear in ancient Egyptian mathematics like on the {\bf Berlin Papyrus 6619}. It contains
the non-primitive triple $6,8,10$ and which is a document dated at about a similar time than the clay tablets. 

\paragraph{}
It has been speculated that experimental exploration of Pythagorean triples also had practical {\bf engineering value} 
because the construction of right angles has architectural or irrigation area measurement applications. 
The interpretation that some of these clay tablets were school tablets indicate that the topic of Pythagorean 
triples must also have been of {\bf educational value}. The Pythagorean triples also paved the way for the 
{\bf Pythagorean theorem}, the statement that $a^2+b^2=c^2$ holds for the sides of a triangle if and only
if the triangle has a right angle. In the remarkable tablet YBC 7289, an isosceles right angle triangle appears with a
rather astounding approximation of $\sqrt{2}$. This was one of the first examples 
for the Pythagorean theorem with non-integer sides but it is also just an example. 

\paragraph{}
Despite many speculations in that direction, there is {\bf no evidence} that the Babylonians were aware of
the Pythagorean theorem. We can {\bf speculate} that they started to {\bf suspect a general rule}. 
We see in the literature and even encyclopedias formulations like 
{\it ``may suggest that the ancient Egyptians knew the Pythagorean theorem"}.
Still, also for such a claim, we lack any historical sources. We have no document on which any such conjecture
is formulated. Formulations like ``may suggest" are a bit reckless as they disregard the difficulty in mathematics
to coming up with conjectures and general rules and then to prove them. There are countless many examples, where 
mathematical rules have been conjectured by looking at small examples and where later, the rule turned out to be false.
Many examples are listed in \cite{stronglawofsmallnumbers}. 
Proto-Pythagorean themes have also appeared also in Chinese documents, including a proof of the Pythagorean 
theorem in the 3-4-5 triangle case which indicates that a general statement was in the air. 
As mathematicians, we know however that stating a fact like $3^2+4^2=5^2$ and visualizing it in a picture is not the 
same than claiming that $a^2+b^2=c^2$ holds in a triangle if one of the angles is a right angle triangle.

\paragraph{}
Jacob Bronowski took in his book \cite{Bronowski} the old fashioned point of view about the discovery of Pythagoras.
He of course was aware about the uncertainty of the sources. But it is the so far best guess that Pythagoras was the
first who proved the theorem. Because of the lack of original documents of Pythagoras, we
might never know who actually proved the Pythagorean theorem the first time. Bronowski tells: 
{\it Pythagoras had thus proved a general theorem: not just for the 3:4:5 triangle of Egypt, or any Babylonian triangle, 
but for every triangle that contains a right angle. He had proved that the square on the longest side or hypotenuse is 
equal to the square on one of the other two sides plus the square on the other if, and only if, the angle they contain 
is a right angle. For instance, the sides 3:4:5 compose a right-angled triangle.
And the same is true of the sides of triangles found by the Babylonians, whether simple as 8:15:17, or forbidding 
as 3367:3456:4825, which leave no doubt that they were good at arithmetic.
To this day, the theorem of Pythagoras remains the most important single theorem in the whole of mathematics. That 
seems a bold and extraordinary thing to say, yet it is not extravagant; because what Pythagoras established is a 
fundamental characterization of the space in which we move, and it is the first time that is translated into numbers. 
And the exact fit of the numbers describes the exact laws that bind the universe.}

\paragraph{}
The ability to make good conjectures and to get a {\bf notion of proof} needed to evolve over time. Already the 
realization that there is a difference between {\bf Examples}, {\bf Conjectures}, {\bf Hypothesis}, {\bf Model} 
and {\bf Theorems} is a cultural achievement. We see the process when watching students learning.
When we learn mathematics we first we confuse the process of proving a result
in general or to just support a phenomenon by giving anecdotal data evidence. 
The ability of {\bf asking questions} like {\bf why} is what 
``makes us curious' \cite{LivioWhy} as it starts a scientific process. 
Because the object {\bf Babylonian graph} seems historically not have appeared in the literature, 
we are in a realm, where we can also observe our own first steps and track misconceptions or mistakes. 

\paragraph{}
Even wrong conjectures can be helpful as they illustrate our status of understanding for a subject. We also are
interested in the Babylonian graph because we are here in a ``data collection" and ``forming conjecture" phase. 
Historians will have to find out whether the Babylonian graph has been mentioned earlier in the literature. 
In a hundred or thousand years, we might know a lot about this graph. Today in 2022, we seem to be in the
``data collection" and ``forming conjectures" phase. In a thousand years, there might be powerful theorems which 
answer all questions. It can of course also be that there is no interest in the object at all building up and that
the topic will remain an obscure example.

\paragraph{}
The origins of the Pythagorean theorem itself is in mystery. Who was the first who conjectured it? This is already
a major step going much beyond just listing examples. The next step is a giant one as it is way beyond conjecture.
Who was the first who proved the Pythagorean theorem? This is difficult as no authentic documents of Pythagoras himself
are known. See \cite{WasPythagorasChinese,Strathern,KahnPythagoras, HiddenHarmonies,PosamentierPythagoras,Maor,Katz2007}. 
Still, if we look at the text evidence, we have to give the Greek mathematicians (and especially Euclid from whom we have
dated documents) the credit to formalize what a ``theorem" and what a ``proof" is and distinguish a ``general statement" 
(which is always true) from a ``statistical statement" (which is true in most cases or with a few exceptions only 
[the nonsensical ``exceptions prove the rule" is even used in colloquial language]) and especially to distinguish 
from ``anecdotal evidence" (which even in our modern times is often mistaken as ``proof" by a mathematically 
untrained person, or then as a crude but effective tool for advertisement or propaganda.) 
The Babylonian triplets were {\bf anecdotal evidence} for the theorem for Pythagoras, not more. It was still
far from a conjecture about a general relation and even further away from a Pythagorean theorem which is a statement
coming with a proof. 

\paragraph{}
There have been a few headlines in the last couple of years claiming that the Babylonians invented trigonomety. 
There is no indication that Babylonians invented trigonometry. This statement depends on 
what ``trigonometry" means. While trigonometry uses ratios of triangles for the definition of the trigonometric
functions, looking at ratios of sides of triangles should not yet count as trigonometry. 
No school curriculum considers that nomenclature when talking about trigonometry. 
Looking at ratios of right angle triangles
is {\bf proto trigonometry at best}. To cite \cite{Robson2002} about research theories in history:
{\it  ``In general, we can say that the successful theory [in the history of mathematics] 
should not only be mathematically valid but historically, archaeologically, and linguistically sensitive too."}
\cite{Mansfield2021} for example has produced the {\bf 2021 controversy}: the paper has been picked up by
media. Math historian Victor Blasj\"o formulated it nicely: `it tricked news outlets into
printing nonsense headlines". 

\section{Experimental explorations}

\paragraph{}
Experimental explorations of a mathematical structure often have predated
theorems considerably. Experiments can lead to {\bf examples which suggest a theorem}.
But there can be a long journey from experiments to theorems. It took a thousand years to get from
{\bf Pythagorean triple explorations} to the {\bf Pythagorean theorem}. Otto Neugebauer already 
speculated that the parametrizations $a=u^2-v^2,b=2uv,c=u^2+v^2$ could have been known thousands of years ago
(\cite{NeugebauerExactSciences} page 39). Bronowski for example gives the example of a pair $(3367,3456)$ \cite{Bronowski} 
which is the parametrization obtained with $u=64,v=27$. The largest number in Plimpton 322 is
$(a,b)=(12709,13500)$ which is obtained with $u=125,v=54$. Obtaining such large numbers without a 
parametrization is harder is not impossible. It needs a bit of patience and some luck. Such examples make it 
likely that the Euclid parametrization was known and used but it is still just a guess.

\paragraph{}
The fact that in the given examples on Plimpton 322, only a few ``random" parameter values $(u,v)$ appear and not
a systematic list, ordered according to $(u,v)$ speak against the knowledge of such a parametrization
but it would be conceivable that some structure was seen like that one number is even and trying $b=2 u v$,
where $u,v$ are factors. 
All primitive triples can be obtained as such and also non-primitives like $6^2+8^2=10^2$ 
have been considered. Finding out what really happened is something for Sherlock Holmes
\cite{Buck1980}, where we see the statement
{\it We can begin by asking if numbers of the form $a^2-b^2$ and $a^2+b^2$ have any special properties.
In doing so, we run the risk of looking at ancient Babylonia from the twentieth century, 
rather than trying to adopt the autochthonous viewpoint.}

\paragraph{}
Creighton Buck further writes:
{\it ``There is no independent information showing that these facts were
known to the Babylonians at the time we conjecture that this tablet was inscribed."}
Indeed, there is no known statement of a result $a^2+b^2=c^2$ for right angle triangles on clay tablets.
Three thousand years ago, it had not been excluded that some
{\bf super large right triangle} with side length $a,b,c$ would satisfy $a^2+b^2 \neq c^2$.
In most tablets, we only see integer side triangles. There is the triple $(1,1,\sqrt{2})$ in YBC 7289
which contains a non-integer side length. We can only imagine how puzzling this must have been. 

\paragraph{}
In the context of finding historical clues, we also can gain insight by looking at what children do. Looking
at early learners is like pointing a telescope to the past. The early steps in mathematics resemble the 
first steps of the pioneers developing the topic. This prompted a historian to claim \cite{Eves}
{\bf "A student should be taught a subject pretty much in the order in which the
subject developed over the ages."} It is a good rule of thumb but of course not universal. Many secrets
from geometry can be appreciated much faster for example when using algebra. The fact that mathematics has
evolved for many thousands of years and in an accelerated way requires a modern student also to pass to the
modern topics faster and taking shortcuts and bypass times of stagnation. 

\paragraph{}
In the context of pedagogy, there is an anecdote of the teacher letting students construct
right angle triangles using paper and ask them measure $a^2+b^2-c^2$. One student group reported in their
presentation that they found a remarkable rule: $a^2+b^2-c^2$ was always small but never zero!
This {\bf anti-Pythagorean theorem} is {\bf academically honest} because every measurement comes
with errors. The students reported {\bf what they measured} and did not report what they {\bf wanted to see}.
It is the most common sin in science to fall into {\bf wishful thinking}. It is a powerful source of 
motivation, but it is dangerous. 
We know that error measurements have a continuous distribution so that without prejudice, it is correct that
in experiments, $a^2+b^2-c^2 \neq 0$ with probability $1$. The students doing the measurements of course had
not been exposed to statistics and data science. A more sophisticated approach would be  to build
a statistical model for the {\bf possible errors}, to make a {\bf hypothesis} and determine the {\bf p-value}, the
evidence against a null hypothesis. A good scientist tries to make the p-value as small as possible and so
give {\bf data evidence} for the Pythagorean theorem. 
This theorem would then be a  {\bf mathematical model}. The scientist then decides whether the measurements 
support the model. This is still far from proving the theorem. To prove the theorem one has 
to placed the statement in a particular frame work, like planar Euclidean geometry. This requires to make
some {\bf idealizations} and {\bf assumptions}. 

\paragraph{}
Also the process of {\bf building a model} or placing a statement in a {\bf particular axiomatic frame work}
is an achievement of Greek mathematics which should not be underestimated. I myself was not taught about
{\bf methods of science} in mathematics but in a philosophy classes. First in high school and later in 
college in a lecture series of Paul Feyerabend. Let me mention the high school part: I have been lucky to have a year of philosophy
in the Schaffhausen highschool with Markus Werner (1944-2016) who was also a successful writer who won a dozen prestigious prizes
like the Herman Hesse literature prize. He started one of the lessons with ``What is the color red?" which led to interesting
discussions about what color is, and whether it is something we can understand. What happens if we mix colors when drawing with 
a yellow and blue crayons what happens if we illuminate an object with blue and yellow lights simultaneously. 
An other of these philosophy lesson started with ``What is the sum of the angles in a triangle?". A student would answer 180 degrees. 
Werner would ask to prove it. An other student would prove it on the board. The class would discuss then what kind of assumptions
went into the proof. Werner would then draw a triangle on a sphere, where the theorem fails. How could we go wrong? 
What was wrong with the proof?  On a sphere, there are triangles where the sum of the three angles is 270 degrees. There 
are 8 triangles on a sphere which partition the sphere up like that. In these 90 degree triangles, the Pythagorean result
$a^2+b^2=c^2$ fails. Actually, for those special 90-90-90 triangles, one has $a^2=b^2=c^2$. They are equilateral 
right angle triangles. Why did the proof which everybody agreed upon fail? What assumptions went into the proof? 

\begin{figure}[!htpb]
\scalebox{0.1}{\includegraphics{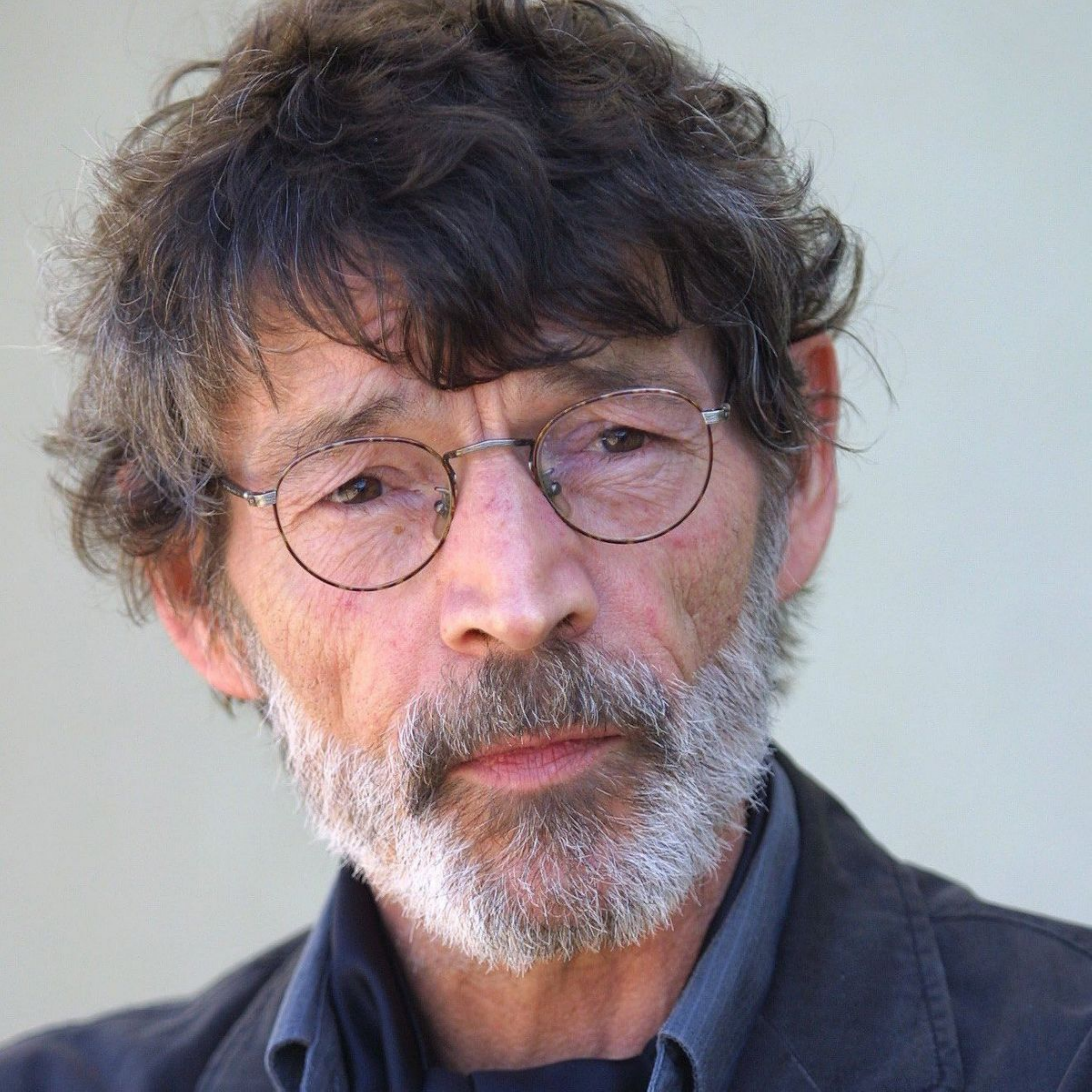}}
\label{adc}
\caption{
Markus Werner (1944-2016), a writer and high school teacher.
}
\end{figure}

\paragraph{}
An even more sophisticated picture appeared since Einstein. Sphere or non-Euclidean geometries are just one of
many Riemannian geometries and Riemannian geometry is a more accurate model of our physical space than Euclidean space. 
\footnote{By the way, Einstein lived 1901/1902 in Schaffhausen for a few months, and lived a few hundred meters from the 
highschool in Schaffhausen, working as a tutor. }
The Pythagorean theorem is well known to fail in our {\bf three dimensional physical space}: it is an
{\bf idealization} dealing with {\bf flat Euclidean space}. 
In a linear algebra setting, assuming a linear flat space, the subject can be dealt with quickly:
define two vectors to be perpendicular if $u \cdot v=0$ and define the length as 
$\sqrt{u \cdot u}$, then check $|u-v|^2=(u-v)\cdot (u-v)=u \cdot u + v \cdot v=|u|^2+|v|^2$.
We make a lot of assumptions although, the result assumes that space is continuous and in particular
that there are perpendicular objects. Then we assume that space has an algebraic structure
in that we can add and scale.

\paragraph{}
But we know since more than 100 years now from general relativity
that every mass bends space and that right angle triangles only satisfy $x^2+y^2=c^2$ in the
{\bf complete absence of matter} or under very special circumstances of the curvature.  But even if we assume
total absence of matter and ignore the presence of virtual particle (which are confirmed 
by phenomena like the Casimir effect), we still do not know because we have no access to any
{\bf Planck scale features} of space. We have no idea what happens if we take a
right angle triangle of side length $a,b,c$ if $a,b,c$ are of the order $10^{-35}$.
Our notions of distance based on measurements using electromagnetic waves do not make sense any more.

\paragraph{}
On a computer small physical distances are no problem- up to some reasonable scale. 
We can for example enjoy looking at features of
the Mandelbrot set on a scale of say $10^{-200}$. It is not difficult for a computer to show us
 topological features of that mathematical object on such a small scale. But we can also with a computer
not explore scales like $10^{-10^{200}}$. 
If the structure of space on the Planck scale would be understood, one can
always ask what happens on an even smaller scale. Once the atom was considered the smallest unit, then 
protons, now we suspect quarks to be part of the smallest ingredients. There was a long way from
speculations by philosophers like Democritus to the current standard model of particle physics. 


\section{About the question}

\paragraph{}
We mentioned the Babylonian graph in our first lecture of Math 22 in the spring of 
2022. It was aimed as an illustration of the fact that 
mathematics is not only {\bf eternal}, but also {\bf infinite}. If you solve one
problem, ten more problems pop up. Having seen in January, the movie ``The eternals" in a 
movie theater, I called the Babylonian graph problem there the {\bf ``Eternals" problem}, because
it had been communicated to us by the eternal Ajak from the Marvel comics universe.
We also used the Babylonian graph as an example in the computer science lecture on May 1st, 2022
Math E 320 to illustrate the process of {\bf experimental mathematics}.
A related graph is the graph in which one takes pairs $a,b$ for which $a+ib$
is a Gaussian prime. We have played with graphs related to number theory also in \cite{Experiments}.

\paragraph{}
In \cite{ExploringCreativity} we looked at the graph $G_n$ with vertex set $V=\{1,2,\dots n\}$, where
two are connected if their sum is a square. So, the connecting rule 
is $a+b=c^2$, (not $a^2+b^2=c^2$). This graph had as a motivation a puzzle
posed by Anna Beliakova of the University of Z\"urich and Dmitrij Nikolenkov of Trogen, a 
high school in Switzerland: {\it Write down the numbers $1-16$ in a row so that the 
sum of two arbitrary neighbors is a square number.} 
This means we have to find a Hamiltonian path in $G_{16}$. Historically, the use of graph theory 
is closely tied to puzzles. William Rowan Hamilton came up with the idea of Hamiltonian paths in the 
context of the {\bf Icosian game}, the problem to find a Hamiltonian cycle on the dodecahedron graph. 
All questions asked for the Babylonian graph can be asked for this {\bf Baliankova-Nikolenkov graph}. 

\paragraph{}
An other class of natural graphs $G_n$ appears on square free integers by connecting two such integers if one divides the other
\cite{PrimesGraphsCohomology}. It has the Euler characteristic $\chi(G_n)) = 1-M(n)$ 
relates to the {\bf Mertens function} $M(x) = \sum_{k=1}^n \mu(k)$ with the M\"obius function $\mu$. 
The value $-\mu(k)$ is the {\bf Poincar\'e-Hopf index} of the vertex $k$ using the Morse function $f(k)=k$.
Adding a new integers is part of a {\bf Morse build-up} because every critical point either has index $1$
or $-1$ and {\bf ``counting" is a Morse theoretical
process}, during which more and more `handles" in the form of topological
balls are added, building an increasingly complex topological structure. I would see later that
this structure has already been studied earlier in \cite{Bjoerner2011}. 

\paragraph{}
In February 2022, we already looked numerically at the 
growth of the graph diameter of the main connected component $B_n'$ of Babylon $B_n$. 
The diameter does not grow monotonically as some parts reconnect. But the experiments suggested
already then that the diameter might increase indefinitely. While easy to see, it  
has not been visible to me at first and I asked it as a question. Only when noticing that there
are connections from integers to a multiple of an integer, the infinite diameter of $B$ became clear
and obvious. In retrospect it now would have been ridiculous to formulate an infinite diameter 
conjecture. 

\paragraph{}
It might well be that one or the other of 
the questions \textcircled{A},\textcircled{B},\textcircled{D}, \textcircled{D} mentioned
here are not difficult to answer. When you look at a new problem the first time, a lot of 
things which later appear ``obvious", are still obscured. It might also be that some
Diophantine problems like the problem of the existence of $K_4$ or $K_5$ subgraphs in $B$
are difficult. The {\bf perfect Euler brick problem} turned out to be hard and evidence that
it is really hard is the fact that it has
remained open for so long. It is well possible that the $K_4$ problem is easy and that
it could be answered by just looking at the problem from the right angle or by having a 
sufficiently strong computer and patience to find one. 

\vfill

\section*{Appendix: A talk on Euler bricks}

\begin{center} \fbox{ \parbox{14cm}{
\begin{center}
{\b Math table, February 24, 2009, Oliver Knill}

\large{ {\bf Treasure Hunting Perfect Euler bricks}} \end{center} }} \end{center}

An Euler brick is a cuboid with integer side dimensions
such that the face diagonals are integers. Already in 1740,
families of Euler bricks have been found. Euler himself constructed
more families. If the space diagonal of an Euler brick is an integer too,
an Euler brick is called a perfect Euler brick. Nobody has found one. There
might be none. Nevertheless, it is an entertaining sport to
go for this treasure hunt for rational cuboids and search - of course with the help of computers.
We especially look in this lecture at the Saunderson parametrization and
give a short proof of a theorem of Spohn \cite{spohn1972} telling that the any of these Euler 
bricks is not perfect. But there are other parameterizations. 

\begin{center}
\scalebox{0.50}{ \includegraphics{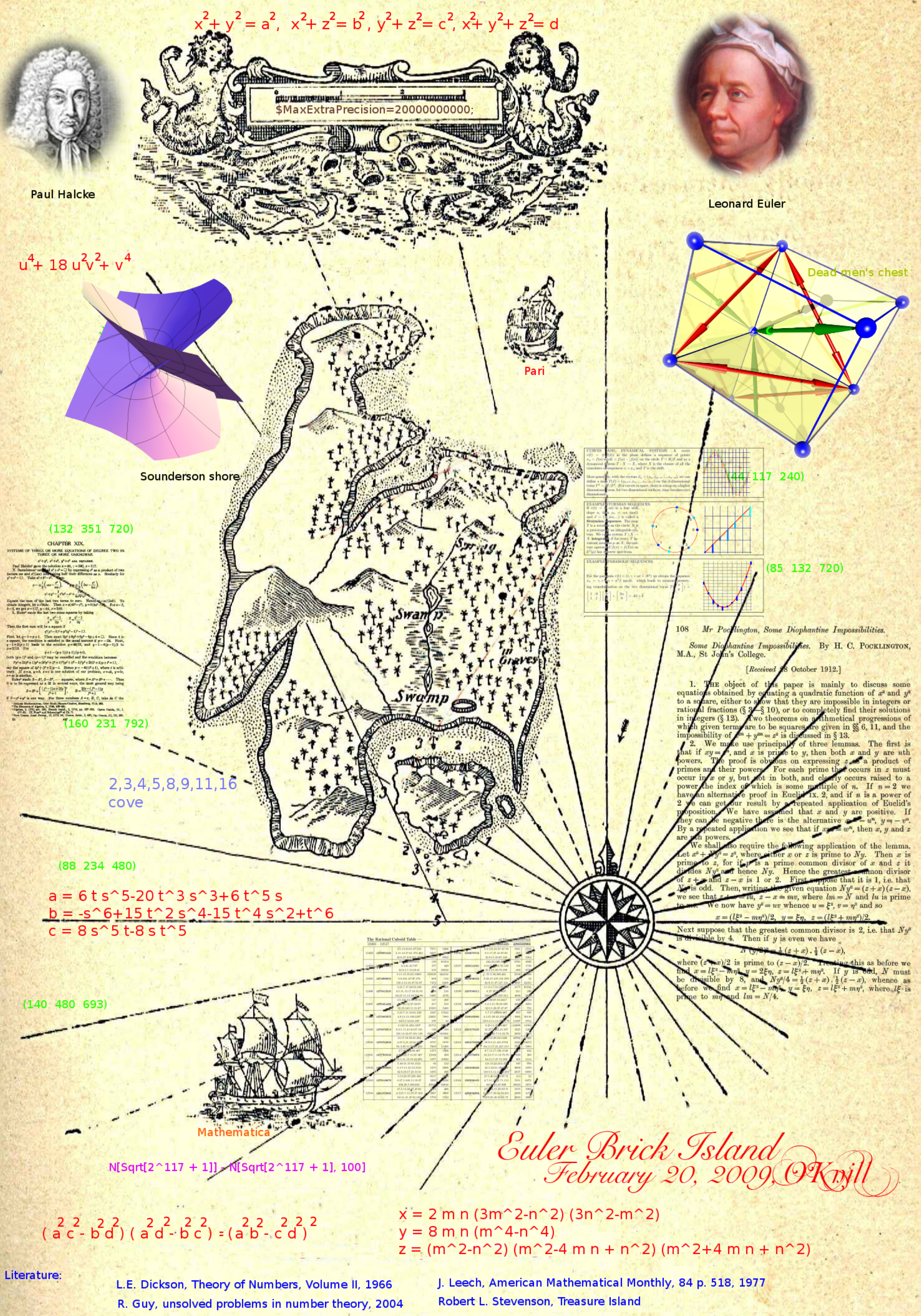}}
\end{center}

\ssection{Introduction: the map of John Flint}

An {\bf Euler brick} is a cuboid of integer side dimensions $a,b,c$ 
such that the face diagonals are integers. If $u,v,w$ are integers satisfying
$u^2+v^2=w^2$, then the Saunderson parametrization
$$ (a,b,c) = ( |u (4v^2-w^2)|, |v (4u^2-w^2)|, |4 u v w|)   $$
leads to an Euler brick.  

\parbox{14cm}{
\begin{center}
\scalebox{0.12}{ \includegraphics{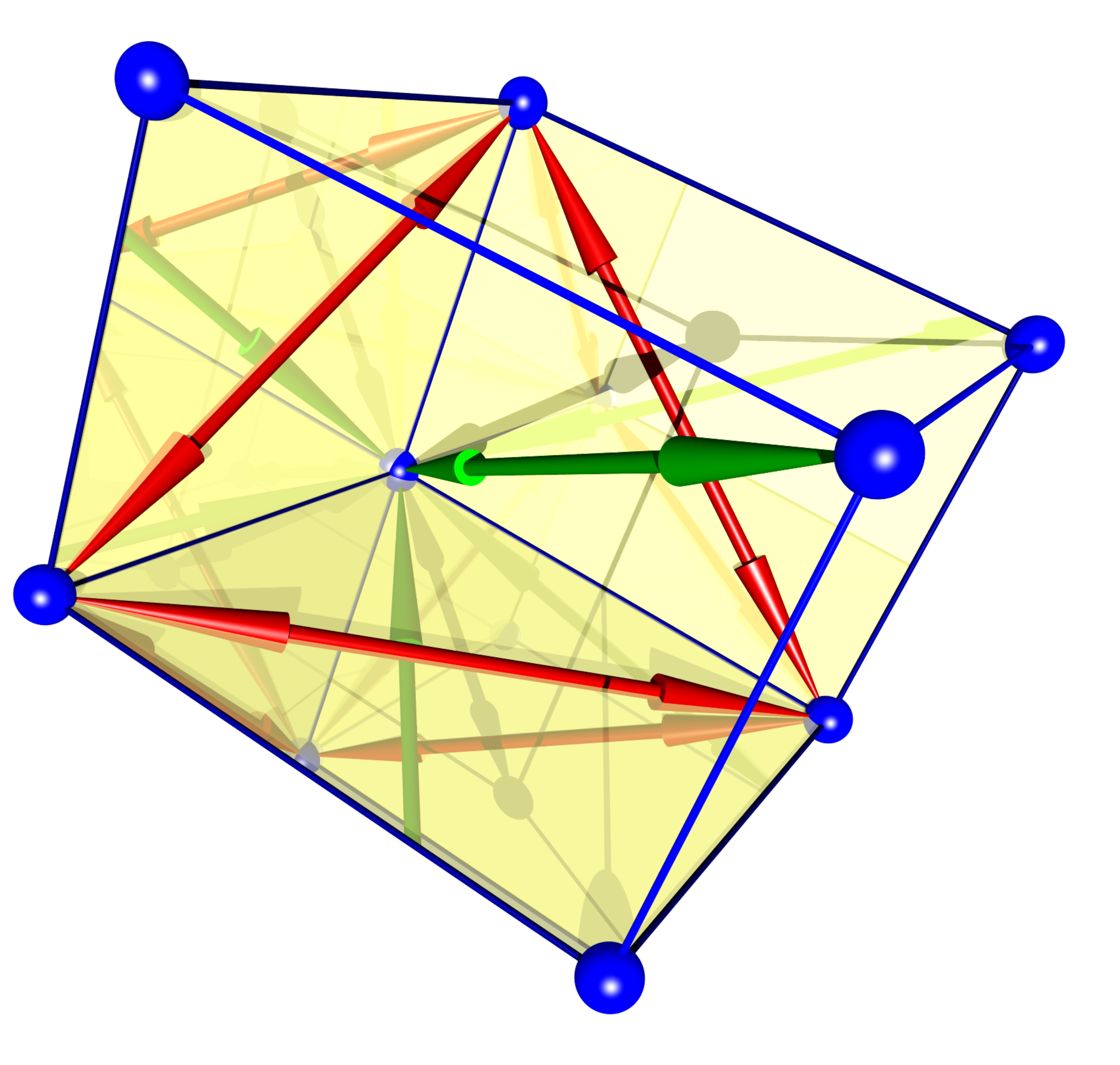}}
\hspace{1cm}
\scalebox{0.50}{ \includegraphics{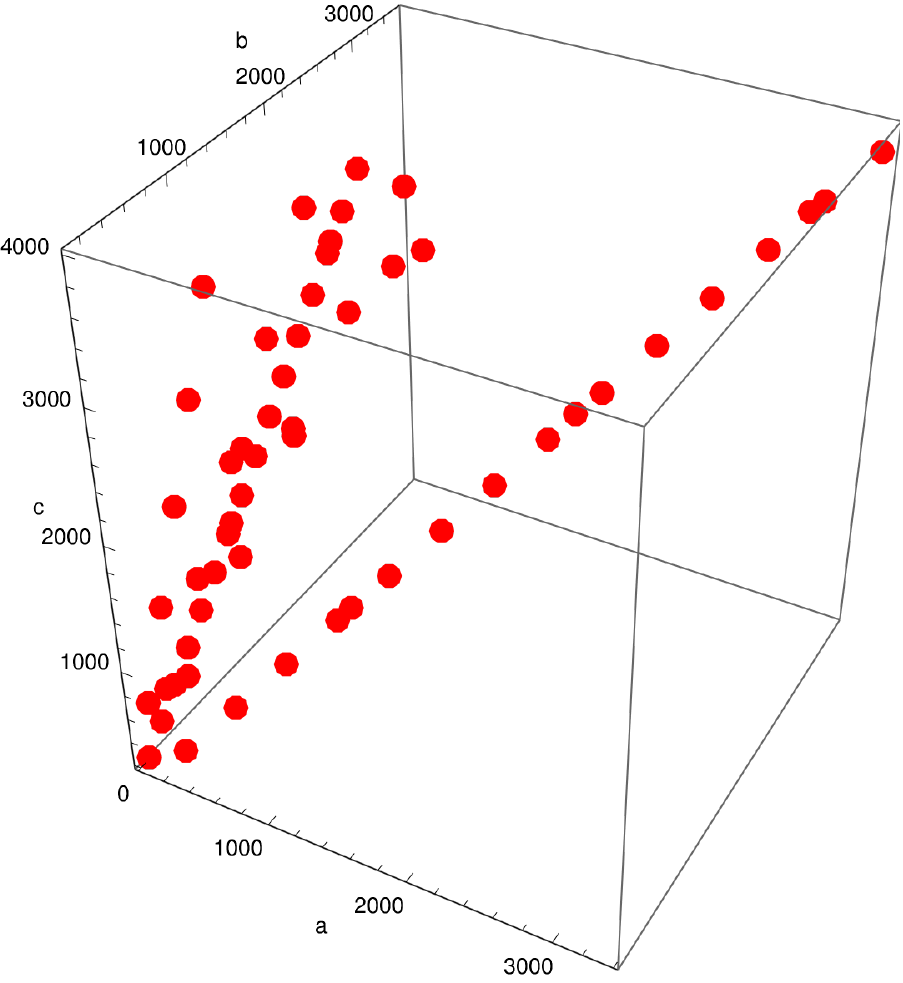}}
\end{center}
}

\parbox{14cm}{
\parbox{6.5cm}{ Fig 1. An Euler brick has integer face diagonals. It is perfect if
                the long diagonal is an integer too. }
\hspace{1cm}
\parbox{6.5cm}{ Fig 2. The smallest Euler bricks $(a,b,c)$ with $a \leq b \leq c$ plotted in the
                parameter space. }
}

\vspace{5mm}

The cuboid with dimensions $(a,b,c) = (240,117,44)$ is an example of an Euler brick. 
It is the smallest Euler brick. It has been found in 1719 by Paul Halcke ( - 1731) \cite{dicksonII}.

If also the space diagonal is an integer, an Euler brick is called
a {\bf perfect Euler brick}. In other words, a cuboid has the properties that the
vertex coordinates and all distances are integers. \\

It is an open mathematical problem, whether a perfect Euler bricks exist.
Nobody has found one, nor proven that it can not exist. One has to 
find integers $(a,b,c)$ such that
$$ \sqrt{a^2+b^2}, \sqrt{a^2+c^2}, \sqrt{b^2+c^2}, \sqrt{a^2+b^2+c^2} $$
are integers. This is called a system of Diophantine equations. 
You can verify yourself that that the Saunderson parametrization produces Euler bricks. 

If we parametrize the Pythagorean triples with
$u= 2s t,  v= s^2-t^2, w = s^2+t^2$, we get
$a = 6 t s^5-20 t^3 s^3+6 t^5 s$,
$b= -s^6+15 t^2 s^4-15 t^4 s^2+t^6$,
$c= 8 s^5 t-8 s t^5$. This defines a parametrized surface
$$ r(s,t) = \langle  6 t s^5-20 t^3 s^3+6 t^5 s,  -s^6+15 t^2 s^4-15 t^4 s^2+t^6,  8 s^5 t-8 s t^5 \rangle $$
which leads for integer $s,t$ to Euler bricks. \\

Indeed, one has then: $a^2+b^2 = (s^2+t^2)^6$, 
$a^2+c^2  = 4 (5s^5t - 6 s^3 t^3 + 5 s t^5)^2$
$b^2+c^2 = (s^6 + 17 s^4 t^2 - 17 s^2 t^4 - t^6)^2$. 

A perfect Euler brick would be obtained if 
$f(t,s) = a^2+b^2+c^2= s^8 + 68*s^6*t^2 - 122*s^4*t^4 + 68*s^2*t^6 + t^8$ were a square.

\pagebreak

\ssection{Brute force search: yo-ho-ho and a bottle of rum!}

There are many Euler bricks which is not parametrized as above:

\parbox{14cm}{
\parbox{6.5cm}{
A brute force search for $1 \leq a,b,c \leq 300$ gives $a=44,b=117,c=240$ 
and $a = 240, b=252 c= 275$ as the only two Euler bricks in that range. 
In the range $1 \leq a<b<c \leq 1000$ there are 10 Euler bricks:  \\

\begin{tabular}{ccc}
    a   &  b  & c   \\ \hline
    44  & 117 & 240 \\
    85  & 132 & 720 \\
    88  & 234 & 480 \\
    132 & 351 & 720 \\
    140 & 480 & 693 \\
    160 & 231 & 792 \\
    176 & 468 & 960 \\
    240 & 252 & 275 \\
    480 & 504 & 550 \\
    720 & 756 & 825 
\end{tabular}
}
\hspace{1.5cm}
\parbox{6.5cm}{
In the $1 \leq a<b<c \leq 2000$,  there are a 15 more, totalling 25.

\begin{tabular}{ccc}
 a  &  b    & c   \\  \hline
170 &  264  &  1440 \\
187 &  1020 &  1584 \\
220 &  585  &  1200 \\
264 &  702  &  1440 \\
280 &  960  &  1386 \\
308 &  819  &  1680 \\
320 &  462  &  1584 \\
352 &  936  &  1920 \\
480 &  504  &  550 \\
720 &  756  &  825 \\
960 &   1008 & 1100 \\
1008 &  1100 &  1155 \\
1200 &  1260 &  1375 \\
1440 &  1512 &  1650 \\
1680 &  1764 &  1925 
\end{tabular}
}
}

Searching $1 \leq a<b<c \leq 4000$, we get $54$ Euler cuboids, in 
$1 \leq a<b<c \leq 8000$ there are 120:

\begin{tiny}
\parbox{15cm}{
\parbox{3.1cm}{
\begin{tabular}{lll|}
44 & 117 & 240\\
85 & 132 & 720\\
88 & 234 & 480\\
132 & 351 & 720\\
140 & 480 & 693\\
160 & 231 & 792\\
170 & 264 & 1440 \\
176 & 468 & 960\\
187 & 1020 & 1584\\
195 & 748 & 6336\\
220 & 585 & 1200\\
240 & 252 & 275\\
255 & 396 & 2160\\
264 & 702 & 1440 \\
280 & 960 & 1386\\
308 & 819 & 1680\\
320 & 462 & 1584\\
340 & 528 & 2880\\
352 & 936 & 1920\\
374 & 2040 & 3168\\
396 & 1053 & 2160 \\
420 & 1440 & 2079\\
425 & 660 & 3600\\
429 & 880 & 2340\\
440 & 1170 & 2400\\
480 & 504 & 550\\
480 & 693 & 2376\\
484 & 1287 & 2640 \\
510 & 792 & 4320\\
528 & 1404 & 2880
\end{tabular} 
}
\parbox{3.1cm}{
\begin{tabular}{lll|} 
528 & 5796 & 6325\\
560 & 1920 & 2772\\
561 & 3060 & 4752\\
572 & 1521 & 3120\\
595 & 924 & 5040 \\
616 & 15 & 3360\\
640 & 924 & 3168\\
660 & 1755 & 3600\\
680 & 1056 & 5760\\
700 & 2400 & 3465\\
704 & 1872 & 3840\\
720 & 756 & 825 \\
748 & 1989 & 4080\\
748 & 4080 & 6336\\
765 & 1188 & 6480\\
780 & 2475 & 2992\\
792 & 2106 & 4320\\
800 & 1155 & 3960\\
828 & 2035 & 3120 \\
832 & 855 & 2640\\
836 & 2223 & 4560\\
840 & 2880 & 4158\\
850 & 1320 & 7200\\
858 & 1760 & 4680\\
880 & 2340 & 4800\\
924 & 2457 & 5040 \\
935 & 1452 & 7920\\
935 & 5100 & 7920\\
960 & 1008 & 1100\\
960 & 1386 & 4752
\end{tabular}
}
\parbox{3.1cm}{
\begin{tabular}{lll|}
968 & 2574 & 5280\\
980 & 3360 & 4851\\
1008 & 1100 & 1155 \\
1012 & 2691 & 5520\\
1056 & 2808 & 5760\\
1100 & 2925 & 6000\\
1120 & 1617 & 5544\\
1120 & 3840 & 5544\\
1144 & 3042 & 6240\\
1155 & 6300 & 6688 \\
1188 & 3159 & 6480\\
1200 & 1260 & 1375\\
1232 & 3276 & 6720\\
1260 & 4320 & 6237\\
1276 & 3393 & 6960\\
1280 & 1848 & 6336\\
1287 & 2640 & 7020 \\
1320 & 3510 & 7200\\
1364 & 3627 & 7440\\
1400 & 4800 & 6930\\
1408 & 3744 & 7680\\
1440 & 1512 & 1650\\
1440 & 2079 & 7128\\
1452 & 3861 & 7920 \\
1540 & 5280 & 7623\\
1560 & 2295 & 5984\\
1560 & 4950 & 5984\\
1600 & 2310 & 7920\\
1656 & 4070 & 6240\\
1664 & 1710 & 5280
\end{tabular}
}
\parbox{3.1cm}{
\begin{tabular}{lll}
1680 & 1764 & 1925 \\
1755 & 4576 & 6732\\
1920 & 2016 & 2200\\
2016 & 2200 & 2310\\
2160 & 2268 & 2475\\
2400 & 2520 & 2750\\
2496 & 2565 & 7920\\
2640 & 2772 & 3025 \\
2880 & 3024 & 3300\\
3024 & 3300 & 3465\\
3120 & 3276 & 3575\\
3360 & 3528 & 3850\\
3600 & 3780 & 4125\\
3840 & 4032 & 4400\\
4032 & 4400 & 4620 \\
4080 & 4284 & 4675\\
4320 & 4536 & 4950\\
4560 & 4788 & 5225\\
4800 & 5040 & 5500\\
5040 & 5292 & 5775\\
5040 & 5500 & 5775\\
5280 & 5544 & 6050 \\
5520 & 5796 & 6325\\
5760 & 6048 & 6600\\
6000 & 6300 & 6875\\
6048 & 6600 & 6930\\
6240 & 6552 & 7150\\
6480 & 6804 & 7425\\
6720 & 7056 & 7700 \\
6960 & 7308 & 7975
\end{tabular}
}
}
\end{tiny}

The number of Euler bricks appears to grow with respect to the box size
because if $(a,b,c)$ is an Euler brick, then $(ka,kb,kc)$ is
an Euler brick too. It would be interesting to know how primitive Euler bricks
are distributed. 


\ssection{Modular considerations: pieces of eight! Pieces of eight! }

If we take a Diophantine equation and consider it modulo some number $n$, then
the equation still holds. Turning things around: if a Diophantine equation has
no solution modulo $n$, then there is no solution in the integers. By checking all
possible solutions in the finite space of all possible cases, we can also determine
some conditions, which have to hold.  \\

Example: $5 x^4 = 3 + 7 y^4$ has no integer solutions because modulo $8$, we have no 
solution because modulo $8$ we have $x^4,y^4 \in \{ 0,1 \}$.  \\

To use this idea, lets assume we deal with prime Euler bricks, bricks for which the greatest
common divisor of $a,b,c$ is $1$.  \\

\ttheorem{
For $n \in \{ 2,3,5,11 \; \}$ as well as $n \in \{ 2^2,3^2,4^2 \; \}$, 
there exists at least one side of an Euler brick which is divisible by $n$. 
}

Proof.  The case $n=2,4,16$ follows directly from properties of Pythagorean triples,
for $n=9$, use that if two (say $x,y$) are divisible by $3$, then $x^2-y^2=a^2-b^2$
is divisible by $9$ and $a=b \; ({\rm mod} \; 3)$ showing that also $z$ has to be divisible by $3$
and the cube is not prime.  \\

\ssection{Searching using irrational rotation: on a dead man's chest}

The problem of solving Diophantine equations has a dynamical system side to it. Take one
of the variables $x$ as time, solve with respect to an other variable say $y$ then 
write $y = (f(x))^{1/n}$ where $f$ is a polynomial.  We can study the dynamical system 
$(f(x))^{1/k} \to f(x+1)^{1/k} \; {\rm mod} \; 1$ and look for $n$ to reach $0$. 
If there are several parameters, we have a dynamical system with multidimensional time.  \\

For the problem to find $s,t$ for which 
$$  \sqrt{s^8+68 t^2 s^6-122 t^4 s^4+68 t^6 s^2+t^8} $$
is close to an integer, we can change the parameter $s,t$ along a line and get incredibly 
close. Unfortunately, we can not hit a lattice point. 

\begin{center} \scalebox{0.80}{ \includegraphics{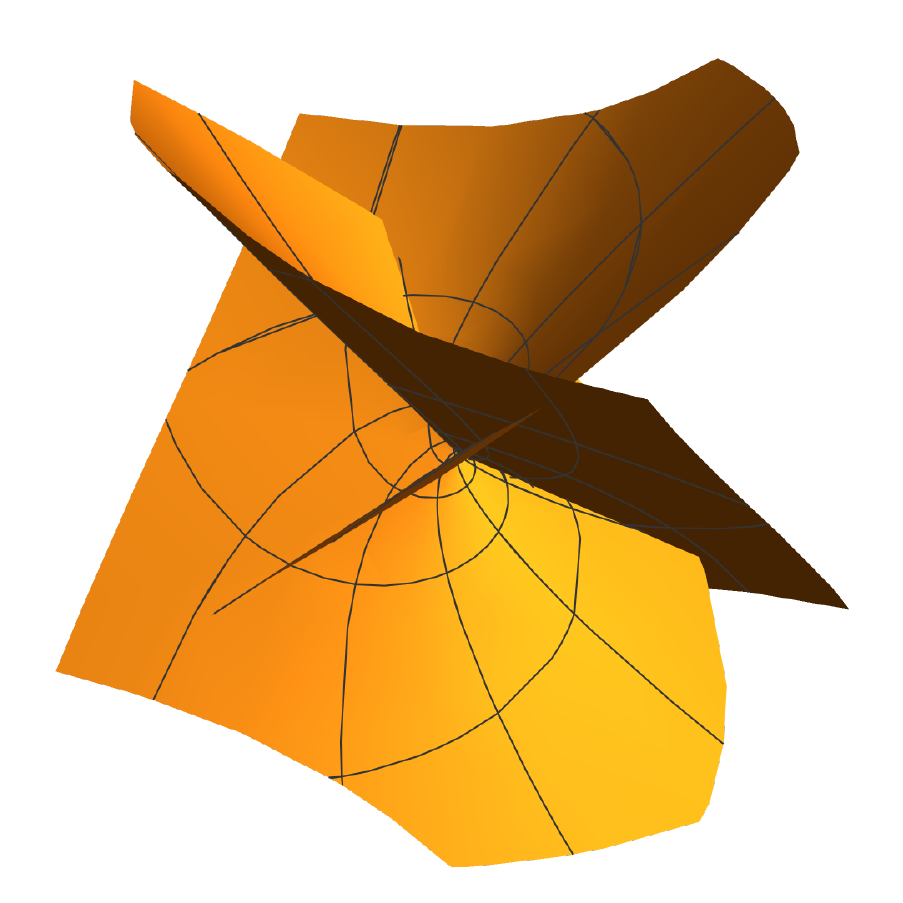} } \end{center}
Fig 1.  The Saunderson surface: a parametrized surface $r(s,t) = (a(s,t),b(s,t),c(s,t))$ of Euler bricks. \\

\ssection{The treasure is not there: ney mate, you are marooned}

Spohn is the "Ben Gunn" of the Euler brick treasure island. He has moved the treasure elsewhere. 
But maybe it does not exist. Anyway, Spohn \cite{spohn1972} proved in 1972: 

\ttheorem{
Theorem (Spohn): There are no perfect Euler bricks on the Saunderson surface of Euler bricks. 
}

Proof. With 
$a = u (4v^2-w^2)$, $b = v (4u^2-w^2)$, 
$c = 4 u v w; w^2=u^2+v^2$, we check
$a^2+b^2+c^2 = w^2 (u^4+18 u^2 v^2 + v^4)$. 
Pocklington \cite{pocklington} has shown first in 1912 that 
$u^4+18 u^2 v^2 + v^4$ can not be be square. His argument is 
more general. We can prove this more easily however: \\

\ttheorem{
Lemma (Pocklington): Unless $x y =0$, the Diophantine equation $x^4+18 x^2 y^2 +y^4 = z^2$ has no solution.
}

Proof: $x,y$ can not have a common factor,  otherwise we could divide it out and 
include it to $z$. Especially, there is no common factor $2$. 
If $x^4+18 x^2 y^2 + y^4  = (x^2+y^2)^2 + 4^2 x^2 y^2  = z^2$ then we have 
Pythagorean triples which can be parametrized.  \\
a) Assume first the triples are primitive, there is no common divisor among the triple
$(x^2+y^2)^2,4^2 x^2 y^2,z^2)$.  \\

(i) If $x,y$ are both odd, we must have
\begin{eqnarray*}
x^2+y^2 &=& 2 u v   \\
4 x y   &=& u^2-v^2  \; . 
\end{eqnarray*}
The first equation proves that $x^2+y^2 =2 \; {\rm mod} \; 4$.
If $2 u v = 2 \; {\rm mod} \; 4$, both $u,v$ must be odd. The second equation can now not be
solved modulo $8$. If $u = 4n \pm 1$, $v = 4m \pm 1$, then $u^2-v^2$ is divisible by $8$.
But the left hand side of the equation is congruent to $4$ modulo $8$.  \\
(ii) If $x$ is odd and $y$ is even, the Pythagorean triple representation is
\begin{eqnarray*}
x^2+y^2 &=& u^2-v^2  \\
2 x y   &=& u v \; . 
\end{eqnarray*}
Because $y$ is even, the second equation shows that $uv$ is divisible by $4$ and because
$u,v$ have no common divisor, wither $u$ is divisible by $4$ or $v$ is divisible by $4$. 
If $u$ is divisible by $4$, the first equation can not be solved modulo $4$. 
If $v$ is divisible by $4$, the first equation has no solution modulo $16$: 
the right hand side is $0,1,4,9$ modulo $16$ 
while the left hand side is congruent to $5,13$ modulo $16$.  \\

b) If there is a common divisor $p$ among $(x^2+y^2)^2$ and $4^2 x^2 y^2$ then it has to be $2$, 
because any other factor $p$ would be a factor
of either $x$ or $y$ as well as of $x^2+y^2$ and so of both $x$ and $y$, which we had excluded
at the very beginning. With a common factor $2$, we have a Pythagorean triple parametrization
\begin{eqnarray*}
x^2+y^2 &=& 4 u v   \\
4 x y   &=& 2(u^2-v^2)  \; . 
\end{eqnarray*}
but since $x,y$ are both odd, $x^2+y^2$ is congruent $2$ modulo $4$ contradicting 
the first equation. 


This finishes the proof of the lemma and so the theorem of Spohn.
It is remarkable that the result of Pocklington does not use infinite decent in this case. 
By the way, the article of Pocklington of 1912 has been checked out many times at Cabot library
since this volume almost falls to dust. 


Side remark: quartic Diophantine equations of this type form an old topic \cite{Mordell1969} (section 4).
Fermat had shown using infinite descent that $u^4+v^4$
is never a square so that $u^4+v^4=z^4$ has no solution. As is well known, he concluded a bit hastily
that he has a proof that $x^p+y^p=z^p$ has no solution for all $p>2$ but that the margin is not large
enough to hold it. 

\ssection{Large numbers: shiver my timbers!}
There are more parametrizations to be explored.  Euler got

\begin{eqnarray*}
a &=& 2 m n (3m^2-n^2) (3n^2-m^2) \\
b &=& 8 m n (m^4-n^4) \\
c &=& (m^2-n^2) (m^2-4 m n + n^2) (m^2+4 m n + n^2)
\end{eqnarray*}
for which $x^2+y^2 = 4m^2n^2(5m^4 - 6m^2n^2 + 5n^4)^2$, 
$^2+z^2 = (m^2+n^2)^6, y^2+z^2 = (m-n)^2 (m+n)^2 (m^4 + 18 m^2 n^2 + n^4)^2$. 

In that case, we have
$x^2+y^2+z^2 = (m^2+n^2)^2 (m^8 + 68 m^6 n^2 - 122 m^4 n^4 + 68 m^2 n^6 + n^8) = 
(m^2+n^2)^2 [ (m^4+n^4)^2  + 2 m^2 n^2 ( 17 m^4 - 31 m^2 n^2 + 17 n^4)]$.  \\

Computer algebra systems like to compute as long as possible in algebraic fields. For example: 

\begin{small} \lstset{language=Mathematica} \lstset{frameround=fttt}
\lstset{backgroundcolor=\color{brightyellow}} \begin{lstlisting}[frame=single]
Expand[(5 + Sqrt[5])^6]
\end{lstlisting} \end{small}

produces the result

\begin{small} \lstset{language=Mathematica} \lstset{frameround=fttt}
\lstset{backgroundcolor=\color{brightyellow}} \begin{lstlisting}[frame=single]
72000+32000 Sqrt[5] 
\end{lstlisting} \end{small}

This is a much more valuable result than a numerical value like
$143554.1753 \dots $. The evaluation of numerical values in Mathematica is
quite mysterious: sometimes, it works quite well:

\begin{small} \lstset{language=Mathematica} \lstset{frameround=fttt}
\lstset{backgroundcolor=\color{brightyellow}} \begin{lstlisting}[frame=single]
N[Sqrt[2^171 + 1]] - N[Sqrt[2^171 + 1], 100]
\end{lstlisting} \end{small}

Sometimes, it does not

\begin{small} \lstset{language=Mathematica} \lstset{frameround=fttt}
\lstset{backgroundcolor=\color{brightyellow}} \begin{lstlisting}[frame=single]
N[Sqrt[2^117 + 1]] - N[Sqrt[2^117 + 1], 100]
\end{lstlisting} \end{small}

which gives in this case a value of $-64$. Even increasing the accuracy
like with 

\begin{small} \lstset{language=Mathematica} \lstset{frameround=fttt}
\lstset{backgroundcolor=\color{brightyellow}} \begin{lstlisting}[frame=single]
$\$$MaxExtraPrecision=20000000000;
\end{lstlisting} \end{small}

Wolfram research promised to fix this problem. \\
By the way, this issue is much better in Pari.
How to compute with large accuracy in the open source algebra system Pari/GP? 
Pari projects algebraic integers correctly, even with millions of digits:

\begin{small} \lstset{language=Mathematica} \lstset{frameround=fttt}
\lstset{backgroundcolor=\color{brightyellow}} \begin{lstlisting}[frame=single]
\p 1000000
a=sqrt(2^117+1) 
\end{lstlisting} \end{small}

It can compute up to 161 million significant digits (you have to increase
the stack size to do so), like defining

\begin{small} \lstset{language=Mathematica} \lstset{frameround=fttt}
\lstset{backgroundcolor=\color{brightyellow}} \begin{lstlisting}[frame=single]
parisize = 800M
\end{lstlisting} \end{small}

in the .gprc file. It still can produce an overflow depending on your machine. 
But working a million digits or so is ok. 

\ssection{History: Captain Flints logbook} 

In 1719 by {\bf Paul Halcke}, a German accountant, who would also do astronomical computations,
found the smallest solution \cite{dicksonII}. Nothing earlier seems to be known. \\

N. Saunderson found in 1740 the parametrization with two parameters mentioned above. 
Only in 1972, it was established by Spohn that the parametrization does not lead to 
that these parametrizations do not lead to perfect Euler bricks. 
Jean Lagrange gave an other argument in 1979 also.  \\

Leonard Euler found in 1770 a second parametrization and in 1772 a third parametrization. 
After his death, more parametrizations were found in his notes. 

\parbox{14cm}{
\begin{center}
\scalebox{0.26}{ \includegraphics{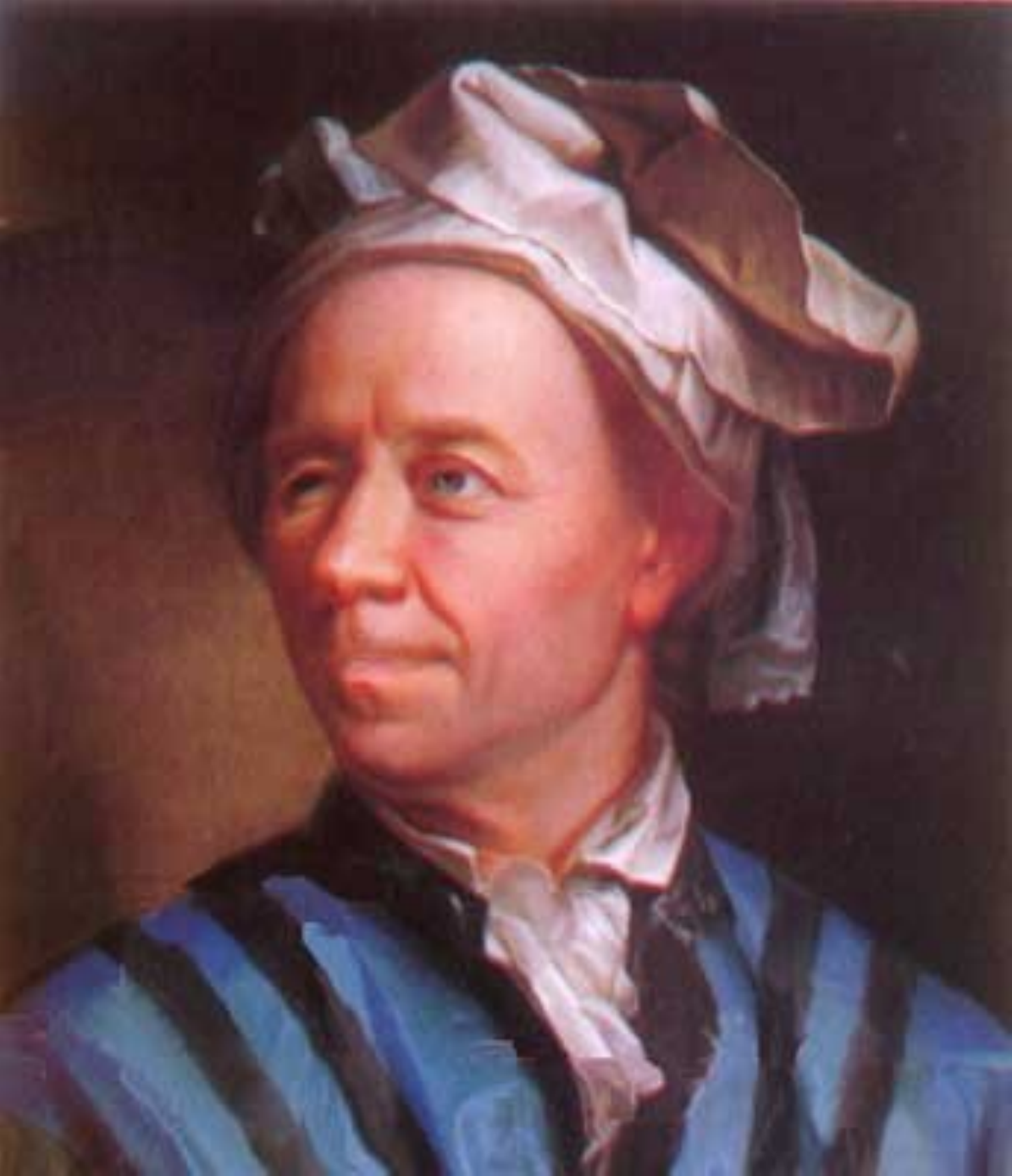}}
\hspace{1cm}
\scalebox{0.49}{ \includegraphics{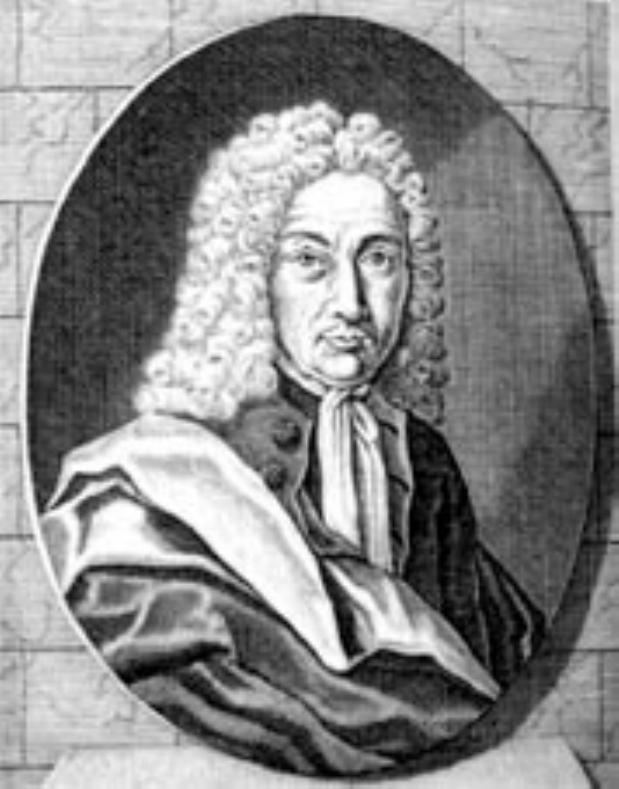}}
\end{center}
}

\ssection{Modern considerations: the black spot}

The topic has appeared several times in American Mathematical Monthly articles and was even a topic for a
PhD theses in 2000 in Europe and 2004 in China. Because of its simplicity, it is certainly of great 
educational value. The topic appears for example in a journal run by undergraduates similar to 
HCMR \cite{ortan}.  \\

Noam Elkies told me: \\

{\it ``The alebraic surface parametrizing Euler bricks is the intersection in $P^5$
of the quadrics $x^2 + y^2 = c^2$, $z^2 + x^2 = b^2$,  $y^2 + z^2 = a^2$
which happens to be a K3 surface of maximal rank, so quite closely
related to much of my own recent work in number theory.  Adding the
condition  $x^2 + y^2 + z^2 = d^2$  yields a surface of general type,
so it might well have no nontrivial rational points but nobody knows
how to prove such a thing."} \\

Noam also remarked that Euler's parametrization would only lead to a finite number of perfect
Cuboids as a consequence of Mordell's theorem. There seems however no reason to be known 
which would tell whether there are maximally finitely many primitive perfect cuboids. \\

There are also relations with elliptic curves since a system of quadratic equations often 
define an elliptic curve. See \cite{macleod}. The article \cite{leechcuboid} which mentions also 
relations with rational points on plane cubic curves.  \\

The problem appeared also in articles for the general public. In 1970  Martin Gardner asked to 
find solutions for which $6$ of the $7$ distances in the cuboid are integers. 
If the large diagonal is an integer, these are no more Euler bricks, unless we would have a perfect
brick.  \\

As for any open problem, it is also interesting to look more fundamental questions.
As with many open problems, the problem to find a perfect Euler brick could be undecidable:
we would not be able to find a proof that there exists no Euler brick. This is
possible only if there is indeed no Euler brick. 
You can read an amusing story about Goldbach conjecture in 
``Uncle Petros and the Goldbach conjecture", where the perspective of such an option
blew all motivation of poor uncle Petros to search for the Goldbach grail. \cite{petros}

\ssection{Treasure problems: scatter and find `em!}

Many unsolved problems like the Goldbach conjecture, the Riemann hypothesis,
the problem to find perfect numbers, or the problem of finding
perfect Euler bricks, finding dense sphere packings in higher dimensions,
are mathematical tasks which could in principle be solved
quickly: by finding an example - if it should exist:

\begin{itemize}
\item Writing down an integer which can not be written as a sum
of two primes would settle the Goldbach conjecture.  

\item Finding an integer for which the sum of the proper divisors is the number itself. 

\item Find a root of the zeta function with $Re(z) \neq 1/2$. Just one lucky punch would be
needed to solve the problem. 
\end{itemize} 

But like treasure hunting, aiming to catch such a treasure is not a good business plan
or a way to make a living: the treasure simply does not need to be there.
If it does not exist, the most skillful treasure hunter can not be successful. \\

But it is the search which is interesting, not the prospect of finding anything. \\

By the way, numerical searches for the grail of a perfect cuboid have been done.
Randal Rathbun has found no perfect cuboid with least edge 
larger than $333750000$. The greatest edge is larger than $10^9$.
See \cite{Guy}. Treasure hunters all over the world have probably gone even further. 
See \cite{rathbun} on ArXiv.  \\

{\bf Update, May 20 2022:} Robert Matson (Matson, Robert D. "Results of a Computer Search for a Perfect Cuboid" (PDF). 
unsolvedproblems.org. Retrieved May 23, 2022.) reports that
{\it there are no perfect cuboids with odd side less than 25 trillion, and no 
perfect cuboids with minimum side less than 500 billion.}


\bibliographystyle{plain}

\end{document}